\documentclass[review]{elsarticle}
\usepackage{lineno,hyperref}
\modulolinenumbers[5]

\usepackage{array}
\usepackage{graphicx}      
\usepackage{caption}        
\usepackage{algorithm}
\usepackage{algorithmic}
\usepackage{subcaption}     
\usepackage{amsmath,amsfonts,amssymb}
\usepackage{epsfig,makecell,float}
\usepackage{setspace,mathrsfs}
\usepackage{tocloft}
\usepackage{textcomp}
\usepackage{multirow,indentfirst,times,color}
\usepackage{booktabs}
\usepackage{threeparttable}
\usepackage{amsthm}
\usepackage{extarrows}

\renewcommand{\arraystretch}{1.5}
\usepackage[titletoc]{appendix}
\usepackage{epstopdf}

\usepackage{natbib}
\biboptions{numbers,sort&compress} 
\usepackage{bm}

\makeatletter
\newif\if@restonecol
\makeatother

\let\clearpage\relax

\journal{Journal of Signal Processing}

\begin{document}\textbf{}

\begin{frontmatter}

\title{
Anti-interrupted sampling repeater jamming via linear canonical Wigner distribution lightweight LFM detection}

\author[mymainaddress,mysecondaryaddress]{Jia-Mian Li}
\author[mymainaddress,mysecondaryaddress]{Bing-Zhao Li\corref{mycorrespondingauthor}}
\cortext[mycorrespondingauthor]{Corresponding author}
\ead{li_bingzhao@bit.edu.cn}

\address[mymainaddress]{School of Mathematics and Statistics, Beijing Institute of Technology, Beijing 100081, China}
\address[mysecondaryaddress]{Beijing Key Laboratory on MCAACI, Beijing Institute of Technology, Beijing 100081, China}

\begin{abstract}
Interrupted sampling repeater jamming (ISRJ) poses a serious threat to radar target detection. Traditional time-frequency (TF) domain anti-jamming methods are prone to TF aliasing in multi-component signal scenarios, and cannot effectively suppress ISRJ with energy close to the real target under low signal-to-noise ratio (SNR) conditions.
To address these challenges, this paper proposes an anti-jamming method based on generalized linear canonical Wigner distribution (GLWD) line detection. By setting the parameters reasonably, the TF image of GLWD can have excellent TF resolution and energy concentration, greatly improving the signal separation and SNR.
Furthermore, in order to enhance the detection capability of the target LFM signal, the existing mobile line segment detection (M-LSD) is improved and the mobile long line segment detection (M-LLSD) is proposed. M-LLSD can detect the target signal more easily and reduce the sensitivity to the jamming signal, so as to efficiently and accurately extract the TF position information of the target signal. Finally, a TF filter is constructed based on the mapping between GLWD and short-time Fourier transform (STFT), performing filtering in the STFT domain to suppress jamming. Simulations and experiments show that the method can effectively suppress such difficult-to-distinguish jamming and is suitable for real-time radar anti-jamming with good robustness.
\end{abstract}

\begin{highlights}
\item The proposed method can suppress ISRJ whose energy and frequency are very close to the real target, which makes up for the shortcomings of the existing anti-ISRJ method.
\item A generalized linear canonical Wigner distribution (GLWD) is proposed, which has superior time-frequency (TF) resolution and noise resistance.
\item The mobile long line segment detection (M-LLSD) is proposed, which can better detect continuous long straight lines and is more suitable for target detection in ISRJ scenarios.
\end{highlights}

\begin{keyword}
Linear canonical transform \sep Fourier transform \sep Linear canonical Wigner distribution \sep Interrupted sampling repeater jamming \sep LFM signal \sep Signal sampling 
\end{keyword}

\end{frontmatter}


\section{Introduction}
Interrupted sampling repeater jamming (ISRJ) based on digital radio frequency memory (DRFM) technology is a coherent deception jamming used in radar electronic warfare \cite{RN28,RN119}. Its principle is that the jammer intercepts the radar transmission signal, samples the signal intermittently, obtains signal slices, and then delays and forwards the slices. ISRJ will produce aliasing with the target echo in the time domain and frequency domain, and will produce a series of false targets after pulse compression (PC), which seriously affects the target detection and identification capabilities of radar. In addition, this method has a fast response speed, can adapt to the rapid changes of radar signals, is easy to deploy on various platforms, and has strong anti-agility capabilities.

Existing anti-ISRJ research has achieved some solutions, which can be divided into: waveform design combined with mismatch filtering \cite{RN120,RN29,RN41}, reconstruction cancellation \cite{RN121,RN32}, time-frequency (TF) domain filtering \cite{RN44,RN33}, deep learning jamming identification \cite{RN34,RN134}. 
The anti-jamming method based on agile waveform and mismatch filter design first optimizes the radar transmission waveform, such as intra-pulse frequency agility and phase encoding to improve the initiative, and then designs the mismatch filter to achieve jamming suppression. The jamming reconstruction elimination method first estimates the jamming sampling width and period, and then iteratively reconstructs and eliminates the jamming signal. Both waveform design and reconstruction cancellation methods rely on the accurate estimation of jamming parameters, which is difficult to cope with the flexibility of DRFM systems. Intelligent methods based on deep learning can automatically mine potential features and are suitable for more complex scenarios. This method has good results and does not require estimation of jamming parameters, but deep learning requires a lot of training costs, which limits the real-time application of the method.
TF domain as one of the important transform domains, it has also received widespread attention in anti-jamming research at home and abroad. The TF domain uses both time domain and frequency domain information, has good signal separation characteristics, and is suitable for practical applications of radar anti-jamming. In \cite{RN129}, an ISRJ suppression method based on energy function detection and bandpass filtering is proposed, which is effective for ISRJ with different parameters. \cite{RN33} The proposed method constructs a max-TF function based on the TF energy distribution of the de-chirped signal, extracts the jamming-free signal segment and generates a filter based on the function.
With the advancement of technology, the jammer can modulate the energy and frequency characteristics of the jamming signal by amplitude frequency and other information, so that the energy and frequency characteristics of the jamming signal are very similar to the target signal \cite{RN147}. The existing TF domain anti-jamming methods rely on the time-domain energy distribution and have obvious defects and limitations. Under the conditions of multi-component signals and low signal-to-noise ratio (SNR), problems such as signal overlap and image blur may occur in the TF plane. Specifically, when the jamming energy is close to or even lower than the target signal energy, the filtering method based on energy distribution may mistakenly classify the target as jamming and filter it out.

Therefore, there is an urgent need for a more accurate radar anti-jamming method that can effectively distinguish between target and jamming and make full use of the differences between target and jamming to achieve better suppression performance.
In response to the above problems, this paper introduces linear canonical transform (LCT) \cite{lct} into the field of anti-jamming, and proposes a new method based on linear canonical Wigner distribution (WD) and lightweight real-time straight line detection. First, a generalized linear canonical Wigner distribution (GLWD) is proposed. This new transformation has four sets of free parameters. By deriving the output SNR model, the parameters are configured to obtain a TF plane with superior TF resolution and noise resistance. Secondly, based on Mobile line segment detection (M-LSD) \cite{gu2021realtime}, this paper proposes Mobile long line segment detection (M-LLSD) by adjusting the loss function and line segment enhancement. This detector can better detect continuous long straight lines and weaken the detection of discontinuous line segments, which is more suitable for target detection in ISRJ scenarios. In addition, the model inherits the efficient single-module framework of M-LSD, has the characteristics of lightweight, and is suitable for practical applications. Considering that WD has no inverse transformation, the target TF distribution position is obtained by M-LLSD detection, and the STFT filter is constructed based on the mapping relationship to filter the echo signal in the STFT domain, and the time domain signal is obtained by inverse transformation to complete jamming suppression.
The main contributions and innovations of this paper are as follows.

(1) A GLWD with excellent time domain resolution and noise resistance is proposed. The superior time domain resolution and noise resistance of GLWD are used to make the overlapping and blurred signals in the traditional TF plane clearly separated in GLWD, which greatly improves the performance upper limit of the TF domain anti-jamming method. In addition, this paper also gives a mathematical proof that the cross term presents discontinuous characteristics in the TF plane.

(2) Based on M-LSD, M-LLSD is proposed by adjusting the loss constraint function. This detection method can better detect continuous long straight lines and weaken the detection of discontinuous line segments, which is more suitable for target detection in ISRJ scenarios. In addition, the model is lightweight and efficient, suitable for practical applications.

(3) Considering that there is no inverse transform for WD, the time domain signal cannot be obtained after filtering for subsequent processing. Therefore, based on the mapping relationship between GLWD and STFT, this paper constructs a TF filter, filters in the STFT domain to suppress jamming, and effectively improves the signal-to-jamming ratio (SJR).


The rest of this paper is arranged as follows: Section \ref{sec2} introduces the signal model and problem formulation. Section \ref{sec3} provides a detailed implementation of the proposed method. Section \ref{sec4} demonstrates the effectiveness of the proposed method in jamming suppression through numerical simulation experiments. Finally, Section \ref{sec5} provides a summary and discussion.

\section{Signal model and problem formulation}
\label{sec2}
\subsection{Signal model}
The radar uses a linear frequency modulation (LFM) signal as the transmitting signal. The LFM signal can be written as
\begin{equation}
s(t) = \text{rect} \left( \frac{t}{T_p} \right) \exp \left( j2 \pi (f_c t+\frac{1}{2}k t^2) \right)
\end{equation}
where $f_c$ is the carrier frequency, $k$ is the frequency modulation slope, $T_p$ is the pulse width, $\text{rect}(\cdot)$ is a rectangular function.
The target echo signal can be expressed as
\begin{equation}
s_{r}(t) = A_{t}\text{rect} \left( \frac{t - \tau_r}{T_p} \right) \exp \left( j2 \pi (f_c (t - \tau_r)+\frac{1}{2}k (t - \tau_r)^2) \right)
\end{equation}
where $A_{t}$ is the amplitude of the target echo signal, $\tau_r$ is the the delay time of target.

Assume that the interrupted sampling signal is a rectangular envelope pulse train, with duration $\tau$ and period $T_s$, The function $p(t)$ is given by \cite{RN119}
\begin{align}
p(t) = \text{rect}\left( \frac{t}{\tau} \right) * \sum_{n=-\infty}^{\infty} \delta(t - nT_s)
\end{align}
\begin{figure}
\centerline{\includegraphics[width=0.4\textwidth]{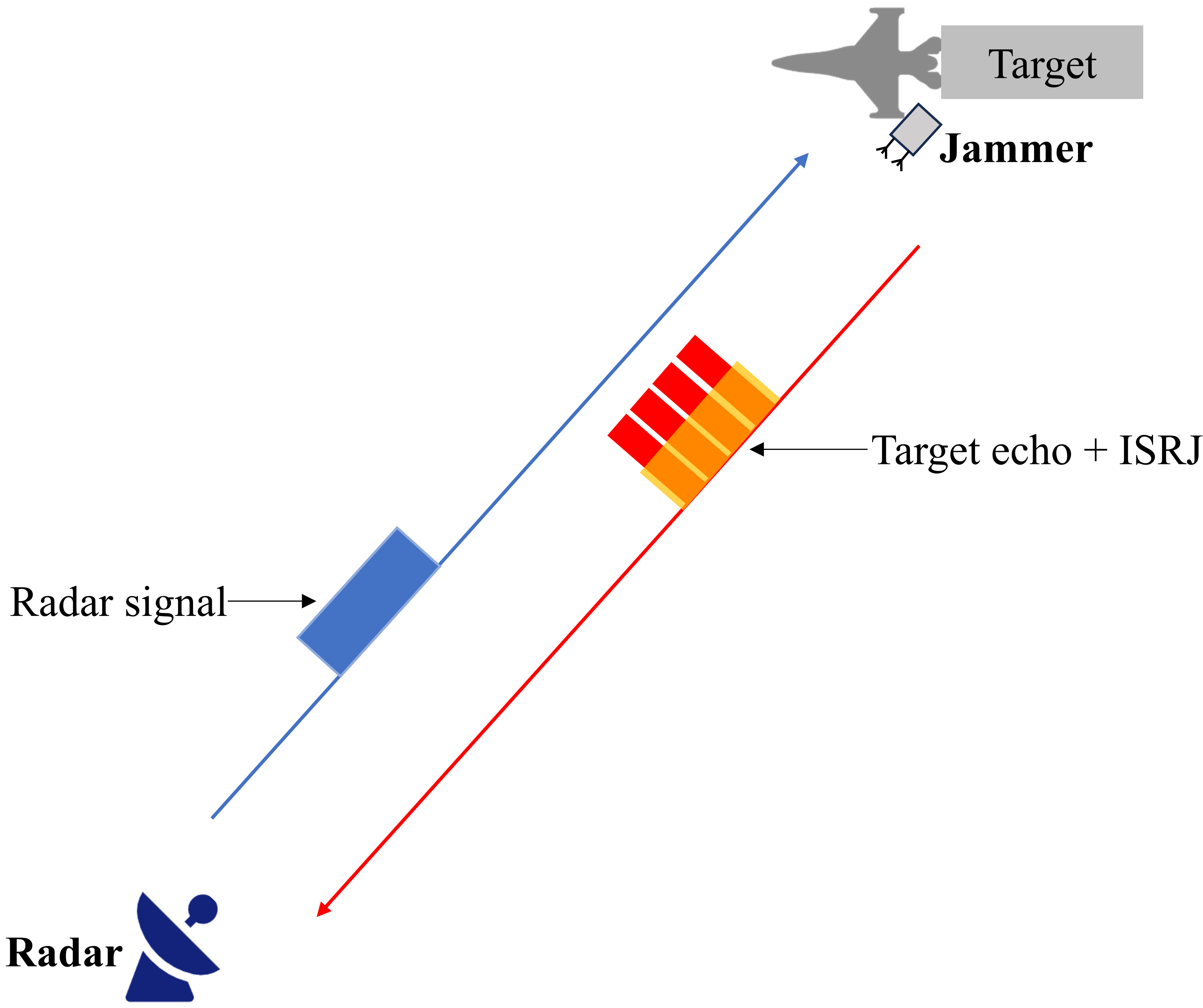}}
\caption{Schematic diagram of radar electronic warfare.}
\end{figure}

The jammer intercepts the radar transmission signal and then multiplies it with the pulse train to obtain a sliced sampling signal
\begin{align}
g(t) = p(t) \cdot s(t)
\end{align}
and its frequency spectrum is
\begin{align}
G(f) = P(f) * S(f) 
=\sum_{n=-\infty}^{+\infty} a_n S(f - n f_s)
\end{align}
where $f_s=1/T_s$, $a_n$ is the scale factor, $P(f)$ and $S(f)$ is the frequency spectrum of $p(t)$ and $s(t)$. 

The jammer intermittently samples the intercepted radar transmission signal and then delays the forwarding to be received by the radar. The jamming signal is
\begin{align}
j(t) = &\sum_{n=1}^{N} \text{rect} \left( \frac{t - \tau_j - (n-1)T_s}{T_j} \right) \notag\\& \times \exp \left( j2 \pi (f_c (t - \tau_j)+\frac{1}{2}k (t - \tau_j)^2) \right)
\end{align}
where $N$ is the number of samples, $T_j$ is the slice width, $T_s$ is the sampling period, $\tau_j$ is the delay time of jammer.

After pulse compression (PC), the jamming signal will form false targets. Since the signal delay and frequency are coupled, the jammer can change the relative distance between the false target and the real target through frequency modulation, so that the false targets will be ahead of or behind the real target. This article will discuss the jamming that is very close to the real target in terms of distance and energy information. This kind of jamming is difficult to suppress by traditional methods. The jamming signal is
\begin{align}
j(t) = A_{j}&\sum_{n=1}^{N} \text{rect} \left( \frac{t - \tau_j - (n-1)T_s}{T_j} \right) \notag\\& \times \exp \left( j2 \pi (f_c (t - \tau_j)+\frac{1}{2}k (t - \tau_j)^2) \right)\exp(j2\pi f_1t)
\end{align}
where $A_j$ is the jamming amplitude, $f_1$ is the frequency shift.
\begin{figure}
\centerline{\includegraphics[width=0.6\textwidth]{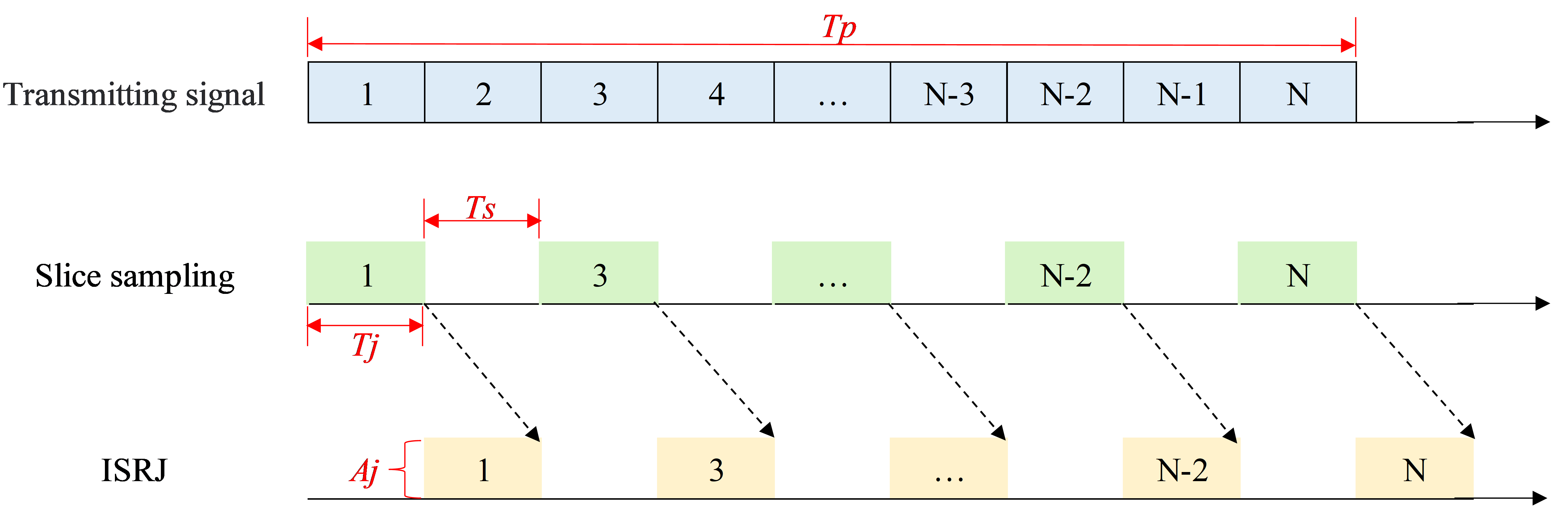}}
\caption{The Principle of ISRJ.}
\end{figure}
The received signal is
\begin{align}
r(t)=s_{r}(t)+j(t)+n(t)
\end{align}
where $n(t)$ is the noise.
\subsection{Wigner distribution}
It is known that the Wigner distributio (WD) play significant role in non-stationary signal processing, especially in detecting LFM signals. \cite{RN138}
To be specific, the WD of \( f \in L^2(\mathbb{R}) \) is defined by \cite{RN125}
\begin{align}
W_f(t, u) = \int_{\mathbb{R}} R_f(t, \tau) e^{-i u \tau} d\tau
\end{align}
where $R_{f,g}(t, \tau) = f(t + \frac{\tau}{2}) g^*(t - \frac{\tau}{2}), R_{f,f}(t, \tau) = R_f(t, \tau)$ and the superscript "*" denotes the complex conjugation.
\subsection{Linear canonical transform}
As a generalization of the Fourier transform (FT) and the fractional Fourier transform (FRFT), the LCT of a signal \( f(t) \) with the parameter matrix \( \mathbf{A} = (a, b; c, d) \)
can be regarded as an affine transformation in the TF plane, which is defined as follows \cite{RN124}
\begin{equation}
F_\mathbf{A}(u) := 
\begin{cases} 
\int_{\mathbb{R}} f(t) K_\mathbf{A}(u, t) dt, & b \neq 0 \\
\sqrt{d} e^{j \frac{cd}{2} u^2} f(d u), & b = 0
\end{cases}
\end{equation}
where
$
K_\mathbf{A}(u, t) := \frac{1}{\sqrt{j2\pi b}} e^{j \left[ \frac{d}{2b} u^2 - \frac{1}{b} ut + \frac{a}{2b} t^2 \right]}.
$
\begin{figure}
\centerline{\includegraphics[width=0.4\textwidth]{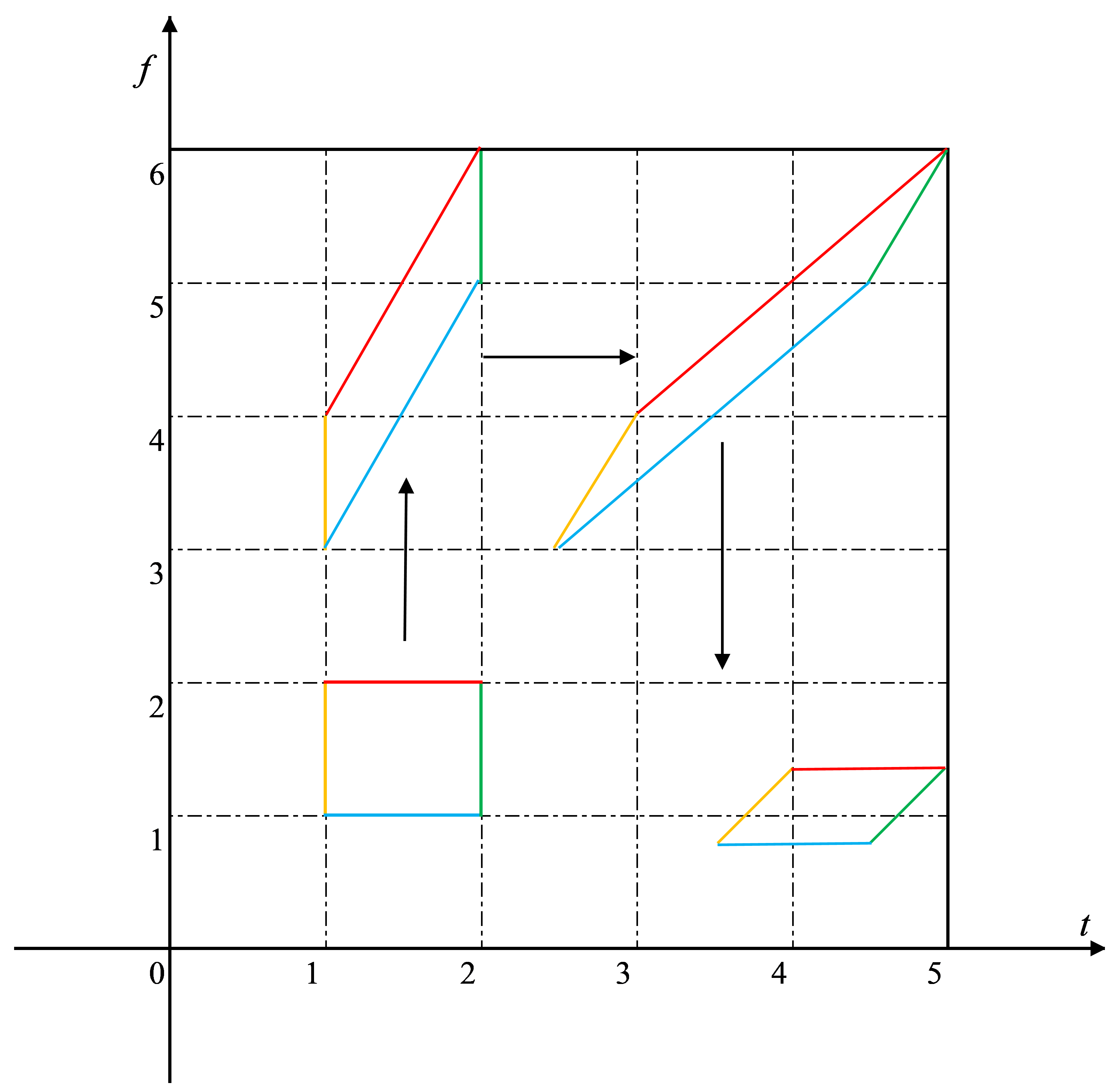}}
\caption{Geometric description of linear canonical transform.}
\label{fig3}
\end{figure}
Compared with the rotation relationship of the signal in the TF plane represented by the FRFT, the LCT represents a more extensive type of affine transformation relationship in the TF plane than the rotation relationship, which can be represented by Fig. ~\ref{fig3}. From Fig. ~\ref{fig3}, it can be seen that under the action of the LCT, the signal not only has a rotation relationship on the TF plane, but also has twisting and stretching changes, but the total support area of the signal on the TF plane will not change. By selecting different LCT parameters $(a, b; c, d)$, the shape and position of the signal on the TF plane can be flexibly changed, which brings greater convenience to the signal processing than the FRFT, so it is more suitable for processing non-stationary signals.
Fast algorithms for LCT \cite{lctfast,RN143,RN144} have also been developed, with the same computational complexity and speed as FFT.

\section{Proposed method}
\label{sec3}
\begin{figure}
\centerline{\includegraphics[width=1\textwidth]{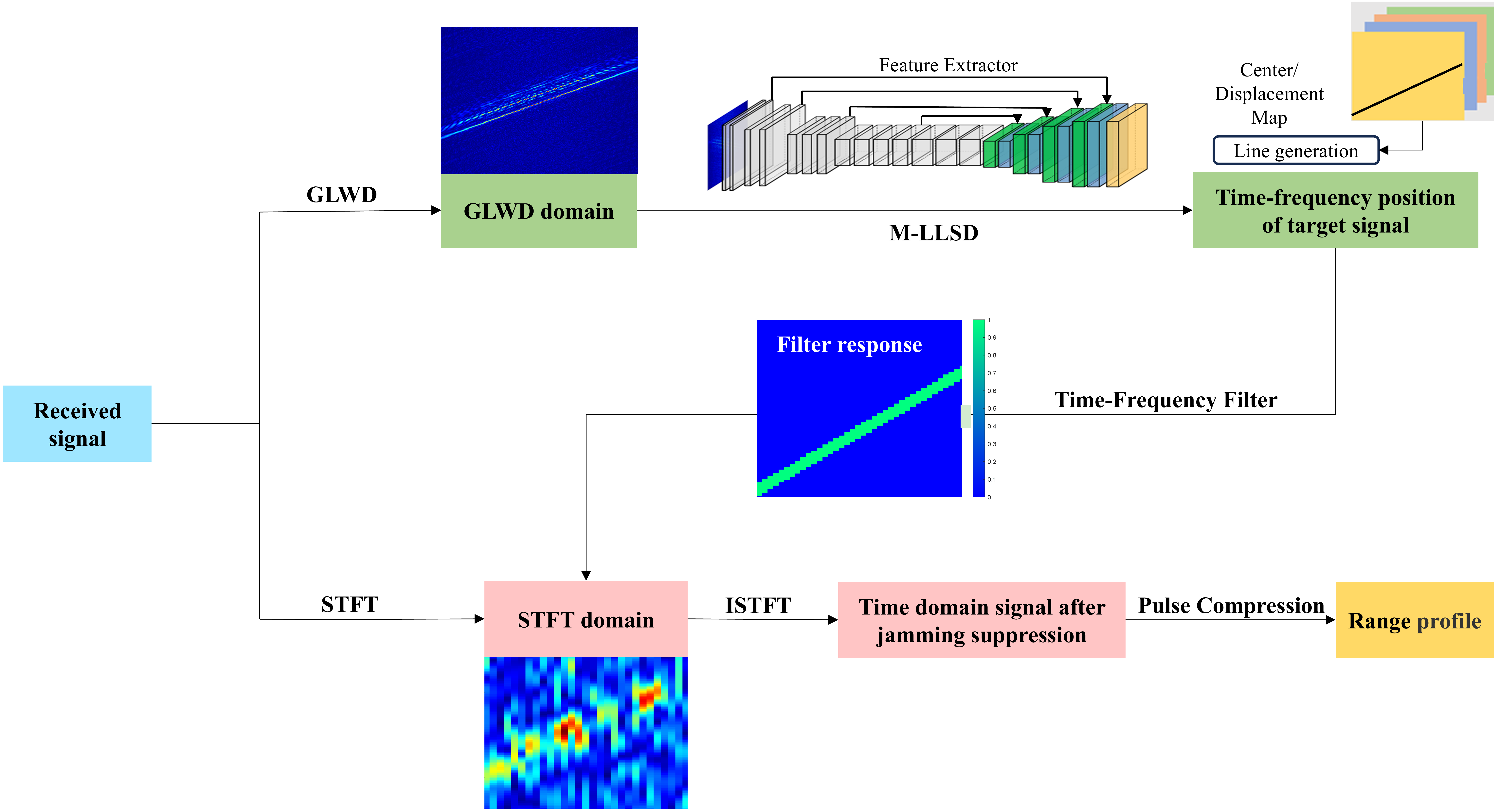}}
\caption{Flowchart of proposed method.}
\label{fig4}
\end{figure}

The main process of the anti-ISRJ method proposed in this paper is shown in Fig. ~\ref{fig4}. First, the TF image is obtained by performing GLWD processing on the received signal. Then, the TF position of the target signal in the received signal is extracted using M-LLSD. In addition, a filter is designed to filter the STFT to eliminate ISRJ, and the time domain signal is obtained by inverse STFT. finally, PC is performed to obtain the radar one-dimensional range profile.
\subsection{Generalized linear canonical Wigner distribution}
For given set of parameters $ \mathbf{A}_i = (a_i, b_i, c_i, d_i), \mathbf{B}_i = (\overline{a}_i,\overline{b}_i, \overline{c}_i, \overline{d}_i) $ satisfy $ a_i d_i - b_i c_i = \overline{a}_i \overline{d}_i - \overline{b}_i \overline{c}_i = 1  (i\in \{1, 2\})$, the new WD associated with LCT of a signal \( f \in L^2(\mathbb{R}) \) are defined as

\begin{align}
W_{f}^{\mathbf{B}_1,\mathbf{B}_2,\mathbf{A}_1,\mathbf{A}_2}(t, u) = \int_{\mathbb{R}} R_{F_{\mathbf{B}_1}, F_{\mathbf{B}_2}}(t, \tau) \times \notag\\& 
\hspace{-4cm} K_{\mathbf{A}_1} \left( u, t + \frac{\tau}{2} \right) K_{\mathbf{A}_2}^* \left(u, t - \frac{\tau}{2} \right) d\tau\label{(11)}
\end{align}
\setlength{\parindent}{15pt}  


Some linear canonical WD methods have been proposed \cite{RN138,RN139,RN140,RN142,RN117}, the newly proposed GLWD can be considered as a generalized form of the existing methods. For example,
when $\mathbf{B}_1$ =  $\mathbf{B}_2$ = $(1, 0; 0, 1)$, $\mathbf{A}_1$ = $\mathbf{A}_2$ = $(0, 1; -1, 0)$, then GLWD is equivalent to WD, when $\mathbf{A}_1$ = $(2a, b; c, \frac{d}{2})$, $\mathbf{A}_2$ = $(-2a, b; c, -\frac{d}{2})$, then GLWD is equivalent to CICFWD \cite{RN117}.
Through flexible parameter configuration, the proposed WD will have stronger noise resistance and higher detection accuracy than existing methods.

This section presents the application of the GLWD in detecting and estimating the parameters of LFM signals. 
Let the LFM signal \( s(t) = e^{i(m t + n t^2)} (m \neq 0) \). 
When the GLWD satisfies the parameter conditions, the LFM signal generates a straight line pulse in the \( (t, u) \) plane.
\begin{equation}
\left\{
\begin{aligned}
    \frac{1}{h_i} &\triangleq 2n \overline{b}_i + \overline{a}_i \neq 0 \quad (i \in \{1, 2\}) \\
    l &\triangleq \frac{\overline{d}_1 - h_1}{8\overline{b}_1} - \frac{\overline{d}_2 - h_2}{8\overline{b}_2} + \frac{a_1 b_2 - a_2 b_1}{8b_1 b_2} = 0
\end{aligned}
\right.
\end{equation}

The amplitude of the GLWD of \( s(t) \) can be presented as follows
\begin{align}
 & \left|{W}_{s}^{\mathbf{B}_{1},\mathbf{B}_{2},\mathbf{A}_{1},\mathbf{A}_{2}}(t,u)\right| \notag\\
 & =\frac{2}{|\frac{1}{b_{1}}+\frac{1}{b_{2}}|}\sqrt{\frac{|h_{1}h_{2}|}{|b_{1}b_{2}|}}\delta\Bigg[u+\frac{1}{\frac{1}{b_{1}}+\frac{1}{b_{2}}}\Bigg(\frac{a_{1}}{b_{1}}+\frac{a_{2}}{b_{2}} \notag\\
 & - \frac{\overline{d}_{1}}{\overline{b}_{1}}-\frac{\overline{d}_{2}}{\overline{b}_{2}}+\frac{h_{1}\overline{b}_{2}+h_{2}\overline{b}_{1}}{\overline{b}_{1}\overline{b}_{2}}\Bigg)t-\frac{h_{1}+h_{2}}{\frac{1}{b_{1}}+\frac{1}{b_{2}}}m\Bigg]
 \label{(13)}
\end{align}
where \( \delta(\cdot) \) denote Dirac delta operator, and the relevant proof is given below.


According to the assumption \( \frac{1}{h_i} \neq 0 \) (\(i \in \{1, 2\}\)) and using Gaussian integral formula $
\int_{R} e^{pt^2 + qt} dt = \sqrt{\frac{\pi}{-p}} e^{-\frac{q^2}{4p}} \quad (p \neq 0, Re(p) \leq 0)
$, we can calculate as follows
\begin{align}
R_{F_{\mathbf{B}_{1}},F_{\mathbf{B}_{2}}}(t,\tau) & =\sqrt{|h_{1}h_{2}|}\mathrm{e}^{i\left(\frac{\overline{d_{1}}}{2b_{1}}-\frac{\overline{d_{2}}}{2b_{2}}\right)t^{2}}\mathrm{e}^{-\frac{ih_{1}\overline{b_{1}}}{2}\left(m-\frac{t}{b_{1}}\right)^{2}}\times \notag\\ 
 & \mathrm{e}^{\frac{ih_{2}\overline{b_{2}}}{2}\left(m-\frac{t}{b_{2}}\right)^{2}}\mathrm{e}^{i\left(\frac{\overline{d_{1}}}{8\overline{b_{1}}}-\frac{\overline{d_{2}}}{8\overline{b_{2}}}-\frac{h_{1}}{8\overline{b_{1}}}+\frac{h_{2}}{8\overline{b_{2}}}\right)\tau^{2}}\times \notag\\
 & \mathrm{e}^{i\left[\frac{\overline{d_{1}}}{2\overline{b_{1}}}t+\frac{\overline{d_{2}}}{2\overline{b_{2}}}t+\frac{h_{1}}{2}\left(m-\frac{t}{\overline{b_{1}}}\right)+\frac{h_{2}}{2}\left(m-\frac{t}{\overline{b_{2}}}\right)\right]\tau}
\end{align}

Therefore, by combining (\ref{(11)}) and (\ref{(14)}), it can be realized that
\begin{align}
\label{(14)}
|W_{s}^{\mathbf{B}_1,\mathbf{B}_2,\mathbf{A}_1,\mathbf{A}_2}(t,u)| &= \frac{1}{2\pi} \sqrt{ \frac{|h_1 h_2|}{|b_1 b_2|}} \Bigg|\int_{\mathbb{R}} e^{-i \left( \frac{1}{2b_1} + \frac{1}{2b_2} \right) u\tau} \,  \nonumber \\
& \hspace{-3cm} \times e^{i l \tau^2} e^{-i \left[ \left( \frac{a_1}{2b_1} + \frac{a_2}{2b_2} - \frac{\overline{d}_1}{2\overline{b}_1} - \frac{\overline{d}_2}{2\overline{b}_2} + \frac{h_1 \overline{b}_2 + h_2 \overline{b}_1}{2\overline{b}_1 \overline{b}_2} \right) t - \frac{(h_1 + h_2)}{2} m \right]\tau} \, d\tau\Bigg|
\end{align}

Obviously, when \( l = 0 \), we arrive at
\begin{align}
\left|W_{s}^{\mathbf{B}_{1}, \mathbf{B}_{2} , \mathbf{A}_{1}, \mathbf{A}_{2}}(t, u)\right| & \nonumber \\ 
& \hspace{-3cm}= \sqrt{\frac{|h_{1} h_{2}|}{|b_{1} b_{2}|}}\delta\left[\left(\frac{1}{2 b_{2}}+\frac{1}{2 b_{1}}\right) u + \left(\frac{a_{1}}{2 b_{1}} + \frac{a_{2}}{2 b_{2}} - \frac{\overline{d_{1}}}{2 \overline{b}_{1}}\right.\right. \nonumber \\
& \hspace{-2.8cm} \left. - \frac{\overline{d_{2}}}{2 \overline{b}_{2}} + \left(\frac{h_{1} \overline b_{2} + h_{2} \overline b_{1}}{2 \overline b_{1} \overline b_{2}}\right) t - \frac{\left(h_{1}+h_{2}\right)}{2} m\right] 
\end{align}
which yields formula (\ref{(13)}). The proof is completed.
\setlength{\parindent}{15pt}  

Let continuous-time signals \( f(t) + n(t) \), where \( n(t) \) represents the noise with zero mean. Then, the mathematical expectation of the GLWD is defined as
\begin{align}
\mathbf{E} \left[ W^{\mathbf{B}_{1}, \mathbf{B}_{2} , \mathbf{A}_{1}, \mathbf{A}_{2}}(t,u) \right] = \nonumber \\
& \hspace{-3cm} W^{\mathbf{B}_{1}, \mathbf{B}_{2} , \mathbf{A}_{1}, \mathbf{A}_{2}}_{f} (t,u) + \mathbf{E} \left[ W^{\mathbf{B}_{1}, \mathbf{B}_{2} , \mathbf{A}_{1}, \mathbf{A}_{2}}_{n}(t,u) \right] 
\end{align}

The output SNR plays an important role in non-stationary signal detection performance improvement. By inheriting the definition of the output SNR of the WD \cite{RN117}, the output SNR of the GLWD can be defined as
\begin{align}
\label{(18)}
\text{OSNR}_\text{GLWD} = \nonumber \notag\\
& \hspace{-2.5cm} \frac{\max_{(t,u) \in \mathbb{R}^2} \left| W_{f}^{\mathbf{B}_{1}, \mathbf{B}_{2} , \mathbf{A}_{1}, \mathbf{A}_{2}}(t,u) \right|}
{\underset{\text{argmax}_{(t,u)}}{\textbf{Mean}} \left| W_{f}^{\mathbf{B}_{1}, \mathbf{B}_{2} , \mathbf{A}_{1}, \mathbf{A}_{2}}(t,u) \right| 
\left\{ \textbf{E} \left[ W_{n}^{\mathbf{B}_{1}, \mathbf{B}_{2} , \mathbf{A}_{1}, \mathbf{A}_{2}}(t,u) \right] \right\}} 
\end{align}
where "Mean" denotes the arithmetic mean for a countable set while the integral average for an uncountable one. 

Assuming that \( D \) is the power spectral density of white noise, we obtain the following result.
The SNR of GLWD for \( f(t) + n(t) \) is
\begin{align}
\label{(19)}
\text{OSNR}_\text{GLWD} = \frac{2\pi }{D} \frac{ |h_1 + h_2 |}{|\frac{1}{b_1} + \frac{1}{b_2}|} 
\end{align}
if and only if \( l = 0 \) and \( \frac{1}{h_i} \neq 0 \) (\( i \in \{1, 2\} \)).

Next, we present the proof of the above formula. Invoking (\ref{(13)}) allows us to conclude that
\begin{align}
\max_{(t,u) \in \mathbb{R}^2} \left| W_{f}^{\mathbf{B}_{1}, \mathbf{B}_{2} , \mathbf{A}_{1}, \mathbf{A}_{2}}(t,u) \right| = \frac{2}{|\frac{1}{b_1} + \frac{1}{b_2}|} \sqrt{\frac{ |h_1 h_2 |}{|{b_1}{b_2}|}}
\end{align}

Meanwhile, the mathematical expectation of zero mean \( n(t) \) can be observed as
\begin{align}
\mathbf{E}\left[ {{W}}_{n}^{\mathbf{B}_{1}, \mathbf{B}_{2} , \mathbf{A}_{1}, \mathbf{A}_{2}}(t,u) \right] &= \notag \\
&\hspace{-0.7cm}\sqrt{\frac{i}{\left( {\overline{a}}_{1}{\overline{b}}_{2}-{\overline{a}}_{2}{\overline{b}}_{1} \right)}} \mathrm{e}^{i{\frac{\zeta_{\mathbf{B}_{1},\mathbf{B}_{2}}+2}{2\left( {\overline{a}}_{1}{\overline{b}}_{2}-{\overline{a}}_{2}{\overline{b}}_{1} \right)}} t^{2}} \notag \\
&\hspace{-1.7cm}\quad \times \frac{D}{\sqrt{2\pi}} \mathcal{K}_{\mathbf{A}_{1}}(u,t) \mathcal{K}_{\mathbf{A}_{2}}^{*}(-u,t) \int_{\mathbb{R}} \mathrm{e}^{i\left[\frac{\zeta_{\mathbf{B}_{1},\mathbf{B}_{2}}-2}{8\left( \overline{a}_{1}\overline{b}_{2}-\overline{a}_{2}\overline{b}_{1} \right)}+\frac{a_{1}}{8b_{1}}-\frac{a_{2}}{8b_{2}}\right] \tau^{2}} \notag \\
&\hspace{-1.7cm}\quad \times \mathrm{e}^{-i\left[ \left( \frac{1}{2b_{1}}+\frac{1}{2b_{2}} \right) u - \left( \frac{\overline{b}_{2}\overline{c}_{1}-\overline{b}_{1}\overline{c}_{2}+\overline{a}_{1}\overline{d}_{2}-\overline{a}_{2}\overline{d}_{1}}{2\left(\overline{a}_{1}\overline{b}_{2}-\overline{a}_{2}\overline{b}_{1}\right)} - \frac{a_{1}}{2b_{1}}-\frac{a_{2}}{2b_{2}} \right) t \right] \tau} \mathrm{d} \tau
\label{(21)}
\end{align}
Therefore
\begin{align}
\left| \mathbf{E} \left[{W}_n^{\mathbf{B}_{1}, \mathbf{B}_{2} , \mathbf{A}_{1}, \mathbf{A}_{2}}(t, u) \right] \right| = \notag \\
&\hspace{-3.7cm}
\frac{D}{\pi} \times
\sqrt{\frac{1}{\left| b_1 b_2 \left( \zeta_{\mathbf{B}_{1},\mathbf{B}_{2}} - 2 \right) + \mu_{\mathbf{A}_{1},\mathbf{A}_{2}} \left( \overline{a}_1 \overline{b}_2 - \overline{a}_2 \overline{b}_1 \right) \right|}}
\label{(22)}
\end{align}
where \( \zeta_{\mathbf{B}_{1},\mathbf{B}_{2}} = \overline{b}_1 \overline{c}_2 + \overline{b}_2 \overline{c}_1 - \overline{a}_1 \overline{d}_2 - \overline{a}_2 \overline{d}_1 \), \( \mu_{\mathbf{A}_{1},\mathbf{A}_{2}} = a_1 b_2 - a_2 b_1 \). Using assumption \( l = 0 \) results in \( \mu_{\mathbf{A}_{1},\mathbf{A}_{2}} = -b_1 b_2 \left( \frac{\overline{d}_1 -h_1}{\overline{b}_1} - \frac{\overline{d}_2- h_2}{\overline{b}_2} \right) \). Then, owing to the above
relation together and conditions \( \overline{a}_i \overline{d}_i = 1 + \overline{b}_i \overline{c}_i , \frac{1}{h_i} = 2n \overline{b}_i + \overline{a}_i({i \in \{1,2\}}) \), we rearrange and simplify the right hand of (\ref{(22)}) to obtain
$
b_1 b_2 \left( \zeta_{\mathbf{B}_{1},\mathbf{B}_{2}} - 2 \right) + \mu_{\mathbf{A}_{1},\mathbf{A}_{2}} \left( \overline{a}_1 \overline{b}_2 - \overline{a}_2 \overline{b}_1 \right) \\= -\frac{b_1 b_2(h_1 + h_2)^2}{h_1 h_2}.
$

Thus, we rewrite (\ref{(22)}) explicitly as
\begin{align}
\left|\mathbf{E} \left[ W^{\mathbf{B}_{1}, \mathbf{B}_{2} , \mathbf{A}_{1}, \mathbf{A}_{2}}(t,u) \right]\right| = \frac{D}{\pi\left|h_{1}+h_{2}\right|} \sqrt{\frac{\left|h_{1} h_{2}\right|}{\left|b_{1} b_{2}\right|}}
\label{(23)}
\end{align}

Finally, substituting (\ref{(21)}) and (\ref{(23)}) into (\ref{(18)}), we obtain (\ref{(19)}).
Next, we propose the expectation-based output SNR inequality model between the GLWD and the WD as well as the GLWD and the CICFWD on a pure deterministic signal added with a zero-mean random noise. \( \text{OSNR}_{\text{WD}} = \frac{2\pi}{D} \), \( \text{OSNR}_{\text{CICFWD}} = \frac{\pi |b|}{D} |h_1 + h_2| \) \cite{RN117}. Obviously, the corresponding conditions for the parameters of GLWD through establishing output SNR inequalities relation between them can be given as follows
\begin{align}
\text{OSNR}_\text{GLWD} > \text{OSNR}_\text{WD} \quad\text{and}\quad
\text{OSNR}_\text{GLWD} > \text{OSNR}_\text{CICFWD} 
\end{align}
Furthermore, it should be noted that the inequalities is equivalent to
\begin{align}
\left| h_1 + h_2 \right| > \left|\frac{1}{b_1} + \frac{1}{b_2}\right| \quad \text{and} \quad \frac{2}{|b|} > \left|\frac{1}{b_1} + \frac{1}{b_2}\right|
\end{align}

As can be seen from the above, GLWD can be regarded as a generalization of WD, mainly reflected in aggregation and anti-noise performance. Through more flexible parameter configuration, better results can be obtained. The target LFM signal and ISRJ are continuous long straight lines and discontinuous line segments on the TF plane, respectively. This characteristic difference lays the foundation for our subsequent target detection. We also need to analyze the cross-term between the two. The cross-term formula is as follows
\begin{align}
W_{s,g}^{\mathbf{B}_1,\mathbf{B}_2,\mathbf{A}_1,\mathbf{A}_2}(t, u) \notag\\& \hspace{-3cm} = \int_{\mathbb{R}} R_{S_{\mathbf{B}_1}, G_{\mathbf{B}_2}}(t, \tau) \times  K_{\mathbf{A}_1} \left( u, t + \frac{\tau}{2} \right) K_{\mathbf{A}_2}^* \left(u, t - \frac{\tau}{2} \right) d\tau
\notag\\& \hspace{-3.0cm} 
=\int_{\mathbb{R}} {S_{\mathbf{B}_1}(t+\frac{\tau}{2}) G^*_{\mathbf{B}_2}}(t-\frac{\tau}{2}) \times  K_{\mathbf{A}_1} \left( u, t + \frac{\tau}{2} \right) K_{\mathbf{A}_2}^* \left(u, t - \frac{\tau}{2} \right) d\tau \notag\\& \hspace{-3cm} 
=\int_{\mathbb{R}} \gamma_1\int_{\mathbb{R}}\text{rect}\left( \frac{t+\frac{\tau}{2}}{T_p} \right)\text{exp}(\alpha_1t^2+\beta_1t)dt \times\gamma_2\int_{\mathbb{R}}\text{rect}\left( \frac{t-\frac{\tau}{2}-nT_s}{T_j} \right)\notag\\& \hspace{-2.5cm}\text{rect}\left( \frac{t-\frac{\tau}{2}}{T_p} \right)\text{exp}(\alpha_2t^2+\beta_2t)dt
\times  K_{\mathbf{A}_1} \left( u, t + \frac{\tau}{2} \right) K_{\mathbf{A}_2}^* \left(u, t - \frac{\tau}{2} \right) d\tau
\label{(26)}
\end{align}
where $\gamma_1$=$\frac{1}{\sqrt{j2{\pi}\overline{b}_{1}}}\text{exp}[j(\frac{k\pi}{4}+\frac{\overline{a}_{1}}{8\overline{b}_{1}})\tau^2+j({\pi}f_c-\frac{u}{2\overline{b}_{1}})\tau+j\frac{\overline{d}_{1}}{2\overline{b}_{1}}u^2]$, $\gamma_2$=$\frac{1}{\sqrt{-j2{\pi}\overline{b}_{2}}}\text{exp}[j(\frac{k\pi}{4}\\+\frac{\overline{a}_{2}}{8\overline{b}_{2}})\tau^2+j(-{\pi}f_c+\frac{u}{2\overline{b}_{2}})\tau+j\frac{\overline{d}_{2}}{2\overline{b}_{2}}u^2]$, $\alpha_1$=$j({\pi}k+\frac{\overline{a_{1}}}{2 \overline{b}_{1}})$, $\alpha_2$=$-j({\pi}k+\frac{\overline{a_{2}}}{2 \overline{b}_{2}})$, $\beta_1$=$j(2{\pi}fc+{\pi}k\tau+\frac{\overline{a_{1}}}{2 \overline{b}_{1}}\tau-\frac{u}{\overline{b}_{1}})$,  $\beta_2$=$j(-2{\pi}fc+{\pi}k\tau+\frac{\overline{a_{2}}}{2 \overline{b}_{2}}\tau+\frac{u}{\overline{b}_{2}})$.

$S_{\mathbf{B}_1}(t+\frac{\tau}{2})$ and $G^*_{\mathbf{B}_2}(t-\frac{\tau}{2})$ in (\ref{(26)}) can be regarded as a $\int_{-\infty}^{+\infty} w({t}) \exp [{j} \theta({t})] {dt}$ integration problem, $w(t)$ is the envelope function and $\theta({t})$ is the signal modulation phase. Using the stationary phase principle, we can solve its approximate solution.

Taking the derivative of the integrand and setting it equal to zero gives
\begin{align}
\frac{d\theta_i({t})}{dt}=2\alpha_it+\beta_i=0,t=-\frac{\beta_i}{2\alpha_i}(i\in \{1, 2\})
\end{align}
\begin{align}
{W}({u})={w}(t=-\frac{\beta_i}{2\alpha_i})
\end{align}

For ISRJ, its envelope function is a rectangular pulse train, and its discontinuous characteristics still exist after integral transformation.
Hence, $R_{S_{\mathbf{B}_1}, G_{\mathbf{B}_2}}(t, \tau)$ still inherits the discontinuous characteristics of ISRJ.
According to the principle of LCT, the integral kernel function can rotate, twist, and stretch the signal on the TF plane, but it will not change the continuity or discontinuity of the signal. Therefore, the energy distribution of the cross-term on the TF plane has an oscillation discontinuity.

\subsection{M-LLSD Model}
\begin{figure}[htbp]
    \centering
    \begin{subfigure}[b]{0.32\textwidth}  
        \centering
        \includegraphics[width=\linewidth]{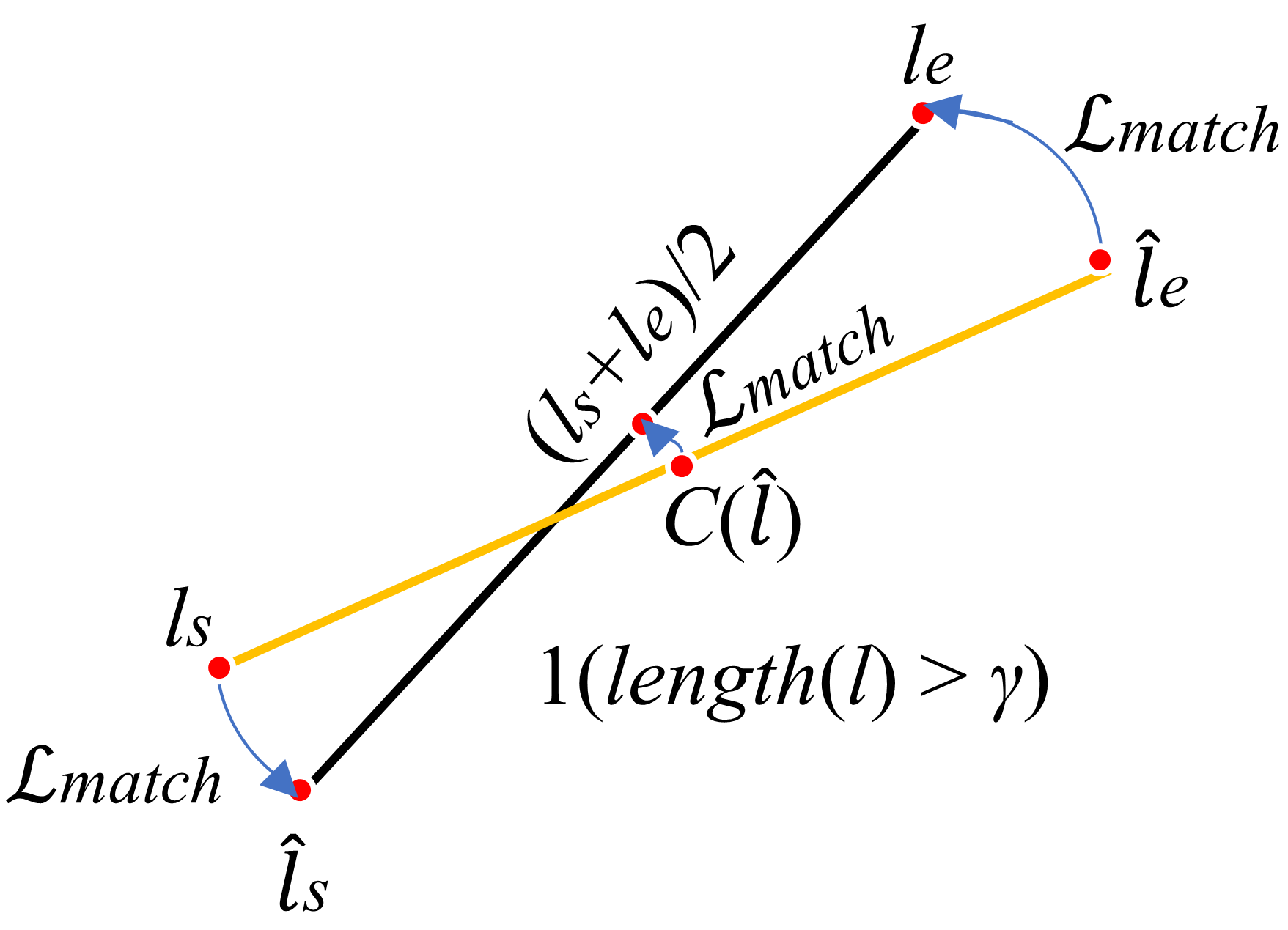}  
        \caption{}
    \end{subfigure}
    \begin{subfigure}[b]{0.25\textwidth}  
        \centering
        \includegraphics[width=\linewidth]{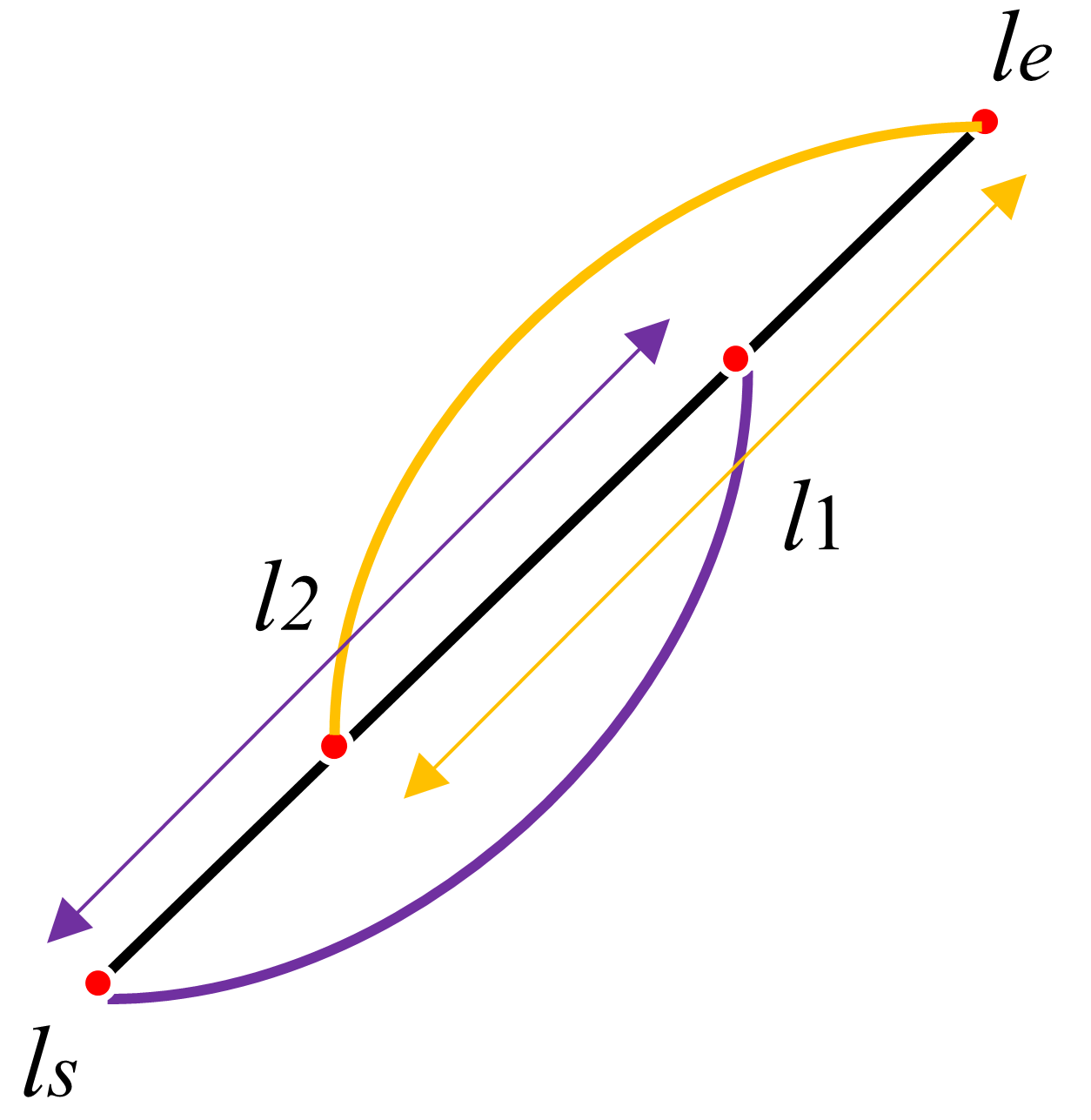}
        \caption{}
    \end{subfigure}
    \begin{subfigure}[b]{0.30\textwidth}  
        \centering
        \includegraphics[width=\linewidth]{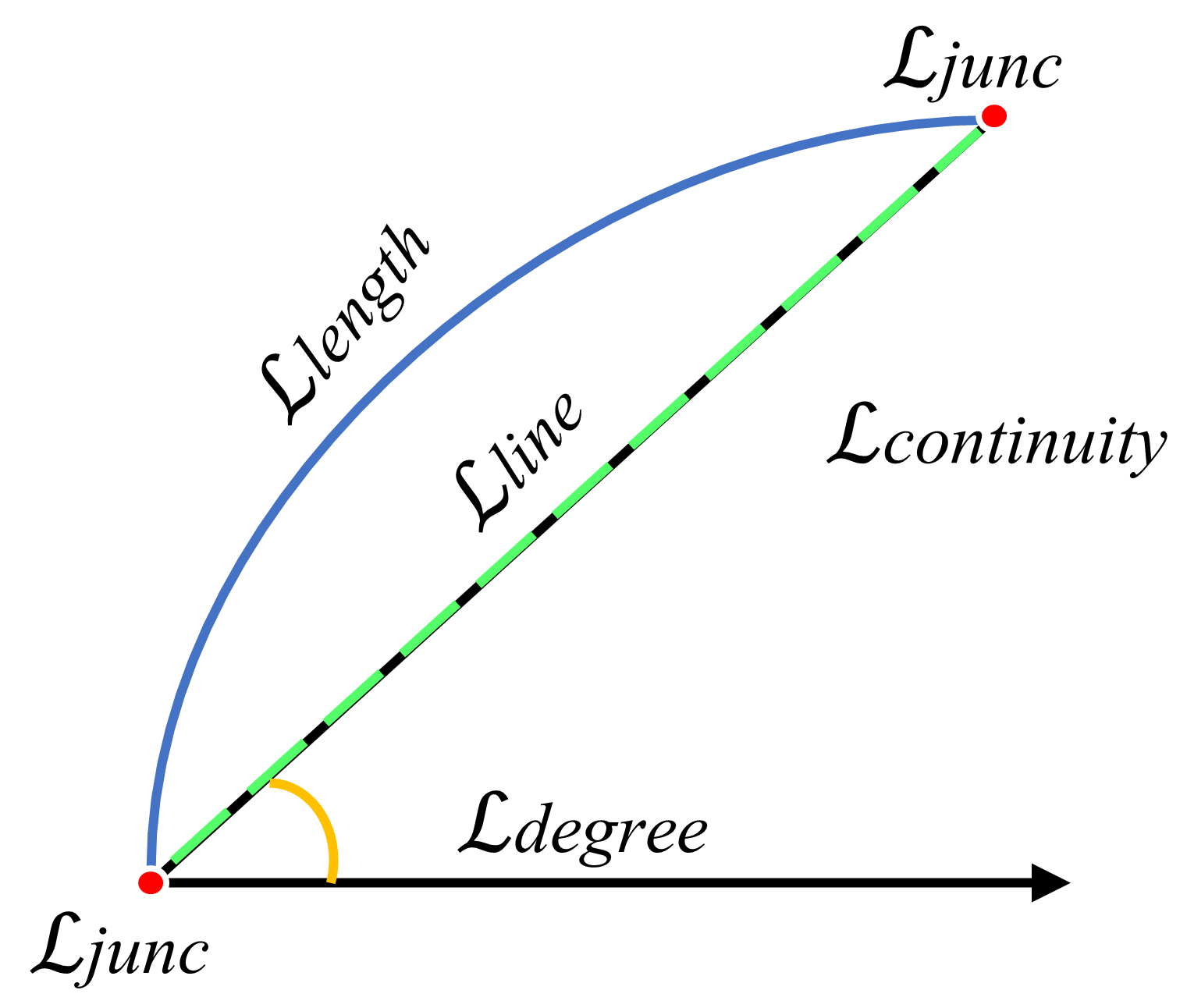}
        \caption{}
    \end{subfigure}
    \caption{Loss function of M-LLSD. (a) Matching loss. (b) SoL augmentation. (c) Geometric loss.}
    \label{fig5}
\end{figure}

M-LSD is a lightweight real-time solution for detecting line segments \cite{RN136,gu2021realtime}. The key idea behind M-LSD is to design an efficient single-module framework, minimizing the backbone network and removing the typical multimodule process for line prediction. To maintain competitive performance in a lightweight network, M-LSD uses a new training scheme: segments of line segment (SoL) augmentation, matching and geometric loss.
Since the target LFM signal appears as a long straight line shape on the GLWD plane, in order to enhance the detection ability of long straight lines, we adjust the SoL augmentation, matching and geometric loss to form M-LLSD, as shown in Fig. ~\ref{fig5}.

(a) Matching Loss: In M-LSD, matching loss is used to optimize the matching degree between predicted line segments and true labels. In order to allow the network to focus on detecting long straight lines and ignore short discontinuous line segments, we introduce a length screening mechanism in the matching loss. Since long straight lines have a larger length than short line segments, a minimum length threshold can be added to the matching loss. Only when the predicted line segment length is greater than this threshold, the loss of the line segment will be considered. If the length of a line segment is less than the threshold, the line segment can be skipped and the loss is not calculated. The modified matching loss $\mathcal{L}_{\text{match}}$ is as follows
\begin{align}
\mathcal{L}_{\text{match}} = & \frac{1}{|\mathbb{M}|} \sum_{(l, \hat{l}) \in \mathbb{M}} \left( \left\|l_{s}-\hat{l}_{s}\right\|_{1} + \left\|l_{e}-\hat{l}_{e}\right\|_{1} \right.\notag \\
& \left. + \left\|\tilde{C}(\hat{l})-\left(l_{s}+l_{e}\right) / 2\right\|_{1} \right)\cdot \mathbf{}{1}(\text{length}(l) > \gamma1)
\end{align}
where $1(\text{length}(l) > \gamma1)
$ is an indicator function that only activates when the length of the segment exceeds a certain threshold $\gamma1$, contributing to the loss calculation. This encourages the model to focus solely on long straight lines.

(b) SoL Augmentation: 
The SoL enhancement method increases the diversity of data by splitting long line segments into multiple sub-segments. In order to prevent the model from detecting intermittent short line segments, the following adjustments can be made modify the splitting rule.
By introducing length threshold $\gamma1$, the line segment is split only when the length is greater than $\gamma1$, the formula is as follows
\begin{align}
k = \left\lfloor \frac{r(l)}{\mu / 2} \right\rfloor - 1 ,\quad \text{if} \quad r(l) > \gamma1
\end{align}

There is an overlap between sub-segments, and the overlap between sub-segments can be increased, which helps the network maintain a strong line segment structure information during training and improves the detection ability of long straight lines. Assuming that the overlap of each sub-segment accounts for a certain proportion of the original line segment, this can be achieved by adjusting the overlap between sub-segments. For example, setting the overlap ratio to $\beta$ (usually between 0.25-0.5), the starting point position of each sub-segment can be calculated as
\begin{align}
l_{s_{i+1}} = l_{e_i} - \beta \cdot (l_{e_i} - l_{s_i}) \quad \text{for each sub-segment}
\end{align}

(c) Geometric Loss:
Add continuity loss $\mathcal{L}_{\text{continuity}}$. $\mathcal{L}_{\text{continuity}}$ is mainly used to judge and optimize the continuity of line segments to ensure that the detected straight lines are continuous in space. Especially during training, we hope to avoid broken line segments or discontinuous parts. If the distance change between two continuous sub-segments exceeds a certain threshold, the error between these sub-segments can be penalized
\begin{align}
\mathcal{L}_{\text{continuity}} = \sum_{\text{l}_i, \text{l}_{i+1}} \mathbf{}{1}(\text{distance}(\text{l}_i, \text{l}_{i+1}) < \gamma2) 
\end{align}
where $\gamma2$ is typically set to a relatively small value.

\subsection{Time-frequency Filtering}
Based on the detection results, it can be determined that the TF interval of the target signal is $(t_1, f_1)-(t_2, f_2)$, which is a straight line with a slope of $k$. Assuming that the width multiple between STFT and GLWD is $w$, the filter can be expressed as follows
\begin{align}
H(t_i, f_i) &=  
\begin{cases}
1, &  f_i \in \left( f_1 + k(t_i - t_1) - \frac{w}{2}, f_1 + k(t_i - t_1) + \frac{w}{2} \right) \\
0, & \text{otherwise}
\end{cases}
\end{align}
The filtered time-domain signal is given by 
\begin{align}
r_{\text{filtered}}(t) = \text{ISTFT}\left(\text{STFT}\left(r(t)\right) \cdot H(t, f)\right)
\label{(34)}
\end{align}


\section{Experimental result and analysis}
\label{sec4}
In this section, the proposed method is evaluated and compared with existing methods, including energy function \cite{RN129} and Max-TF \cite{RN33}.
The radar electronic warfare scenario simulated in this paper is: the radar transmits LFM waveform, and the target platform is equipped with a DRFM jammer. The jammer can intercept the radar transmission signal and generate jamming signals through modulation to be received by the radar, affecting the target detection, and the noise power in the environment is large.
In addition, the jammer can make the energy and frequency of the jamming false target very similar to the real target through frequency shift and amplitude modulation, thereby increasing the difficulty of jamming suppression. This chapter simulates a single ISRJ scenario and a compound ISRJ scenario, and verifies the effectiveness of the proposed method through experiments.
This paper uses signal-to-jamming ratio improvement factor (SJRIF) \cite{RN130}, SJR, and the signal loss ratio (SLR) \cite{RN132} related indicators, defined as
\begin{align}
\text{SJR}=20\mathrm{log}_{10}\left(\frac{{A}_s}{{A}_j}\right),\text{SLR}=20\mathrm{log}_{10}\left(\frac{\overline{A}_s}{{A}_s}\right)
\end{align}
\begin{align}
\text{SJRIF}=20\mathrm{log}_{10}\left(\frac{\overline{A}_s}{\overline{A}_j}\right)-20\mathrm{log}_{10}\left(\frac{A_s}{A_j}\right)
\end{align}
where ${A}_s$ and $\overline{A}_s$ are the amplitudes of the real target before and after jamming suppression, respectively, and ${A}_j$ and $\overline{A}_j$ are the amplitudes of the maximum amplitude false target before and after jamming suppression, respectively.
\begin{table}[ht]
\small 
    \begin{center} 
        \begin{minipage}{0.7\textwidth} 
            \centering
            \caption{Simulation Parameters.} 
            \label{table1}
            \small 
            \renewcommand{\arraystretch}{1}  
            \begin{tabular}{@{}p{0.6\textwidth}cl@{}}
                \toprule
                \text{Parameter} & \text{Value} \\ 
                \midrule
                Waveform pulse duration & 1 $\mu$s \\
                Waveform bandwidth & 100 MHz \\
                Carrier frequency & 10 GHz \\
                Pulse repetition frequency & 5 KHz \\
                Targets' location & 900 m \\
                Main false targets' location & 927 m \\
                ISRJ sampling period & 0.1875$\mu$s\\
                ISRJ slice width & 0.125$\mu$s \\
                SJR & 0 dB \\
                SNR & -12 dB \\
                \bottomrule
            \end{tabular}
        \end{minipage}
    \end{center}
\end{table}
\subsection{Single ISRJ}
\label{sec:single_isrj}
\begin{figure}[htbp]
    \centering
    \begin{subfigure}[b]{0.3\textwidth}  
        \centering
        \includegraphics[width=\linewidth]{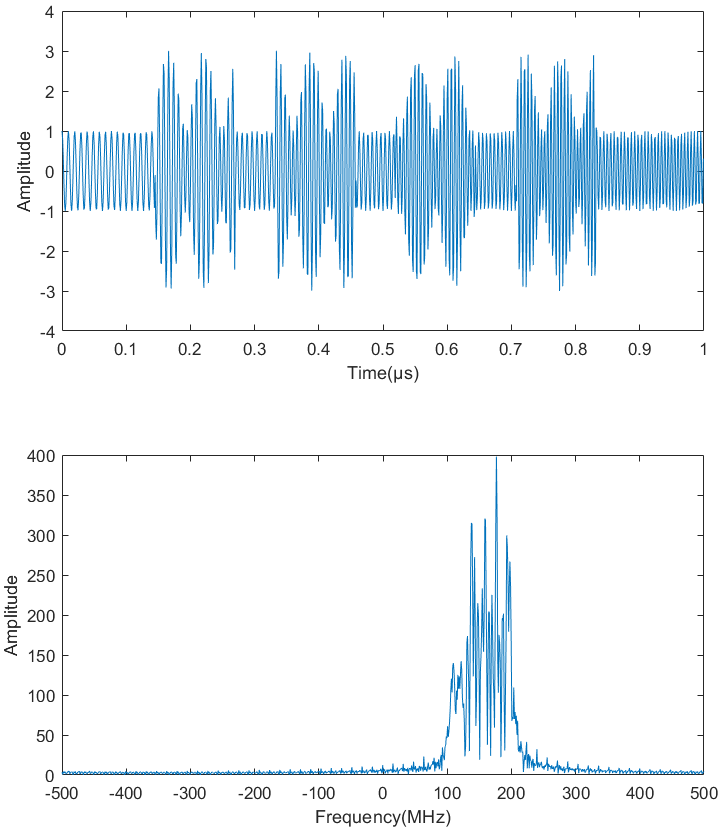}  
        \caption{}
    \end{subfigure}
    \begin{subfigure}[b]{0.3\textwidth}  
        \centering
        \includegraphics[width=\linewidth]{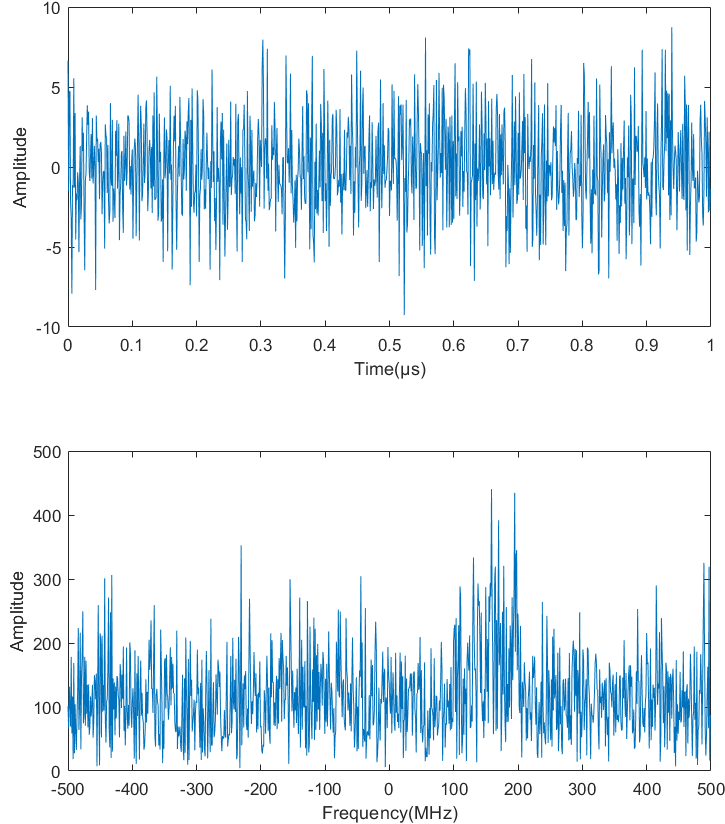}
        \caption{}
    \end{subfigure}
       \caption{Time and frequency domains of signals in a single ISRJ scenario. (a) Target+ISRJ. (b) Target+ISRJ+Noise.}
    \label{fig6}
\end{figure}
\begin{figure}[htbp]
    \centering
    \begin{subfigure}[b]{0.3\textwidth}  
        \centering
        \includegraphics[width=\linewidth]{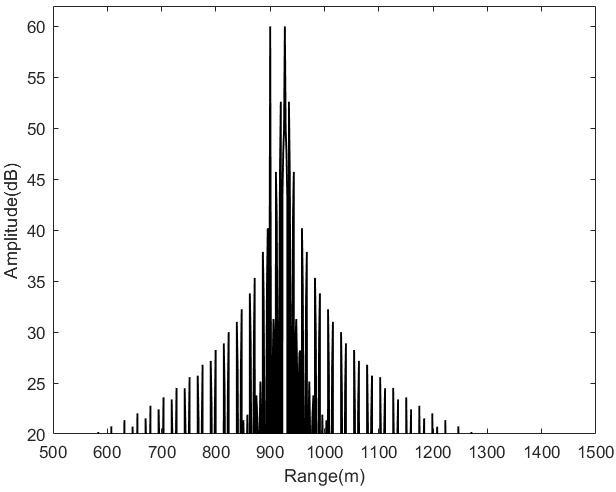}  
        \caption{}
    \end{subfigure}
    \begin{subfigure}[b]{0.3\textwidth}  
        \centering
        \includegraphics[width=\linewidth]{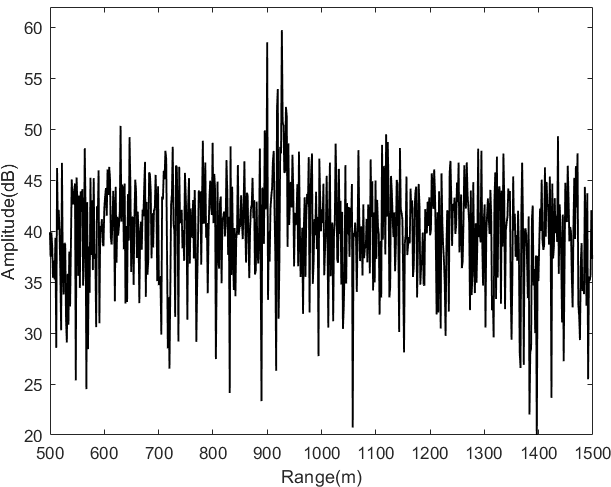}
        \caption{}
    \end{subfigure}
       \caption{Pulse compression in a single ISRJ scenario. (a) Target+ISRJ. (b) Target+ISRJ+Noise.}
    \label{fig7}
\end{figure}
\begin{figure}[htbp]
    \centering
    \begin{subfigure}[b]{0.3\textwidth}  
        \centering
        \includegraphics[width=\linewidth]{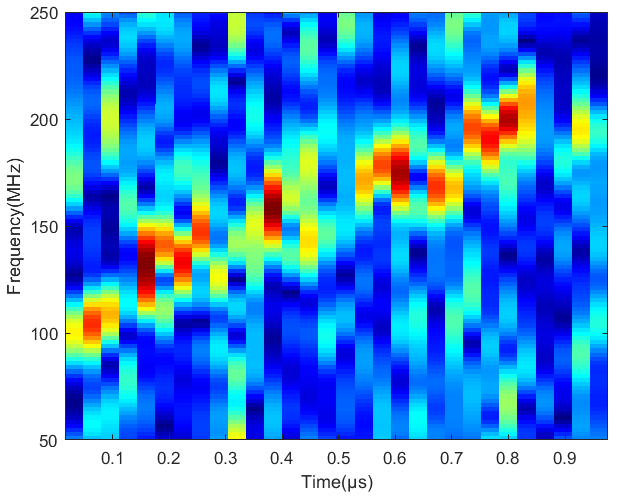}  
        \caption{}
    \end{subfigure}
    \begin{subfigure}[b]{0.3\textwidth}  
        \centering
        \includegraphics[width=\linewidth]{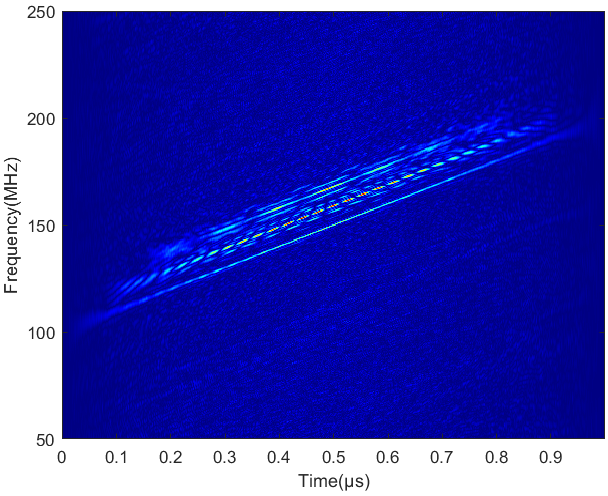}
        \caption{}
    \end{subfigure}
       \caption{Time-frequency image. (a) STFT. (b) GLWD.}
    \label{fig8}
\end{figure}

The relevant parameter settings are shown in Table \ref{table1}. Most existing methods set the jamming energy significantly higher than the target, while this paper sets the jamming energy close to the target energy. As shown in Fig. ~\ref{fig6}, it can be seen that after adding noise, it is difficult to distinguish the ISRJ segment and the target segment in the time domain signal, which makes the method \cite{RN129} relying on the time domain energy function invalid. As shown in Fig. ~\ref{fig7}(a), the energy of the main jamming target is equal to the target, and the range position is very close to the target. As shown in Fig. ~\ref{fig7}(b), after adding noise, the energy of the jamming target is slightly higher than the real target, which increases the difficulty of distinction.
As shown in Fig. ~\ref{fig8}(a), in this case, the TF image obtained after STFT of the received signal is very blurred, and the position of the target signal cannot be distinguished from it. The TF image obtained after GLWD processing, as shown in Fig. ~\ref{fig8}(b), shows that the noise is effectively suppressed, and the target, jamming and cross-term parts are clearly visible.

\begin{figure}
\centerline{\includegraphics[width=0.6\textwidth]{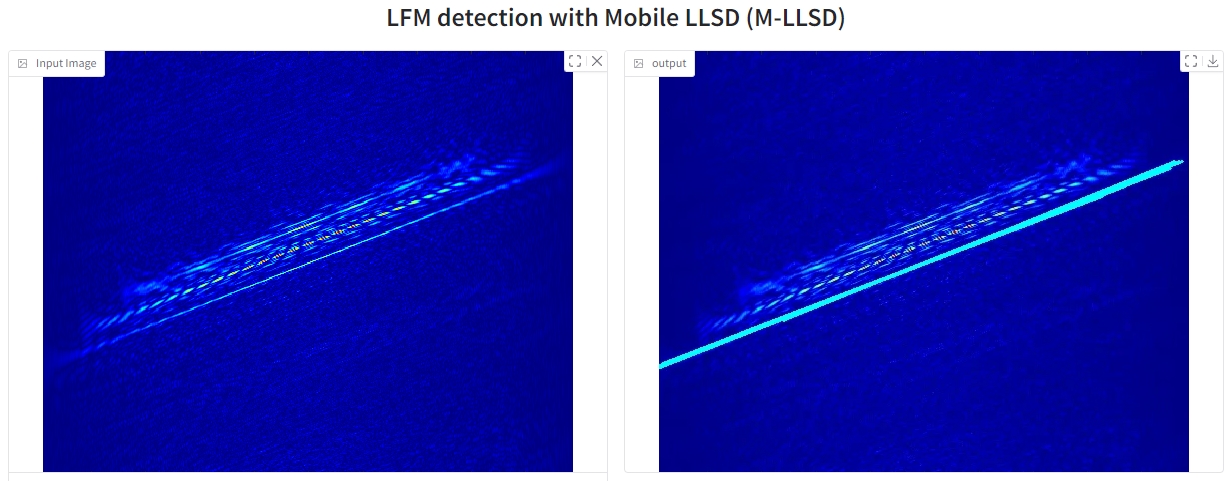}}
\caption{Line detection of target signal.}
\label{fig9}
\end{figure}

\begin{figure}[htbp]
    \centering
    \begin{subfigure}[b]{0.32\textwidth}  
        \centering
        \includegraphics[width=\linewidth]{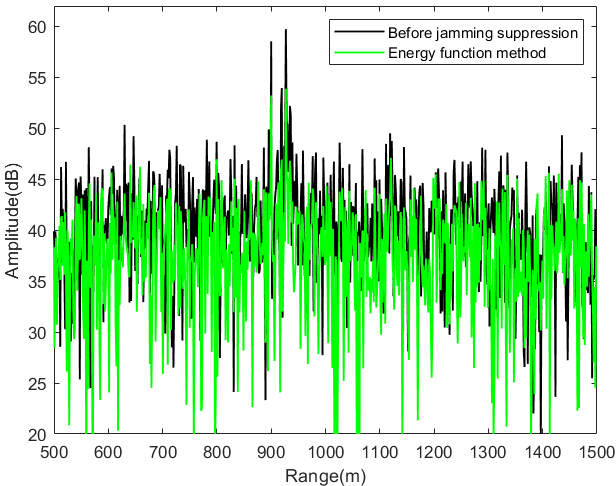}  
        \caption{}
    \end{subfigure}
    \begin{subfigure}[b]{0.32\textwidth}  
        \centering
        \includegraphics[width=\linewidth]{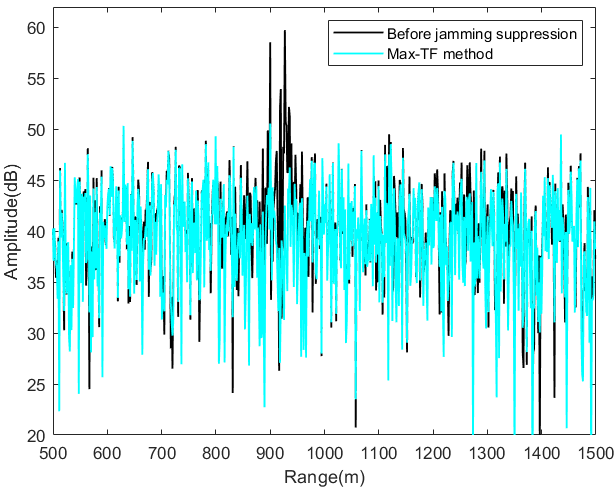}
        \caption{}
    \end{subfigure}
    \begin{subfigure}[b]{0.32\textwidth}  
        \centering
        \includegraphics[width=\linewidth]{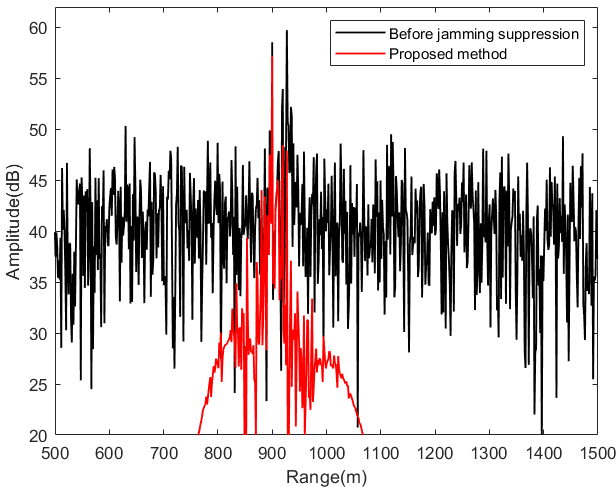}
        \caption{}
    \end{subfigure}
    \caption{Jamming suppression results. (a) Energy function. (b) Max-TF. (c) Proposed method.}
    \label{fig10}
\end{figure}

\begin{figure}[htbp]
    \centering
    \begin{subfigure}[b]{0.24\textwidth}  
        \centering
        \includegraphics[width=\linewidth]{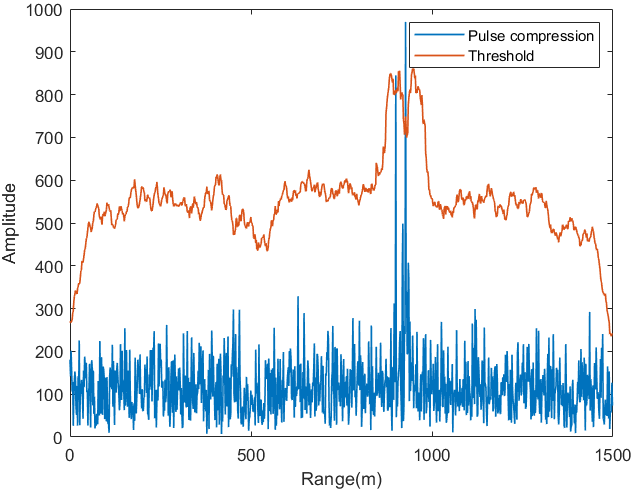}  
        \caption{}
    \end{subfigure}
    \begin{subfigure}[b]{0.24\textwidth}  
        \centering
        \includegraphics[width=\linewidth]{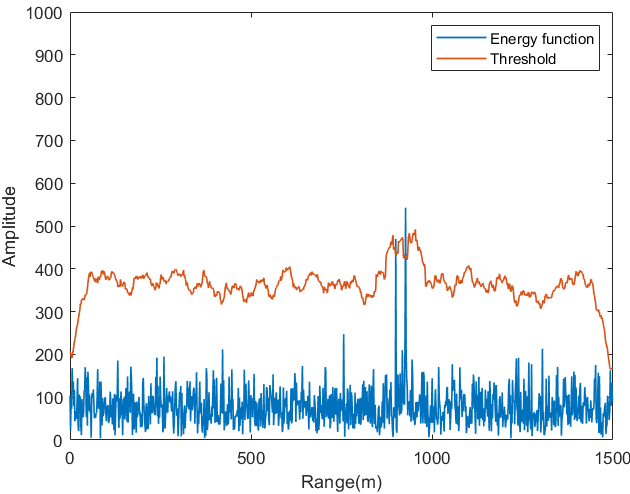}
        \caption{}
    \end{subfigure}
    \begin{subfigure}[b]{0.24\textwidth}  
        \centering
        \includegraphics[width=\linewidth]{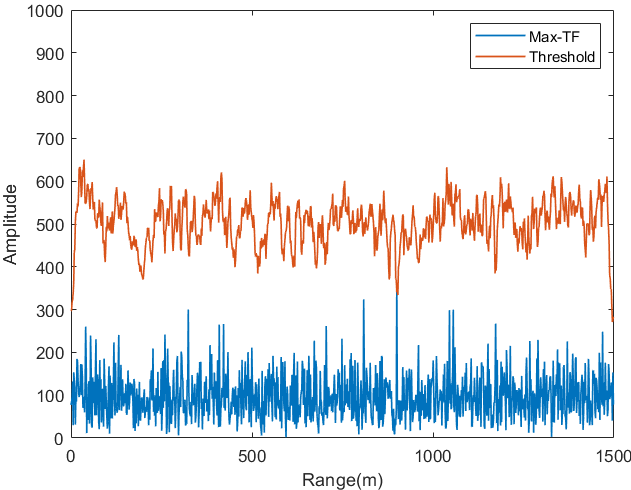}
        \caption{}
    \end{subfigure}
    \begin{subfigure}[b]{0.24\textwidth}  
        \centering
        \includegraphics[width=\linewidth]{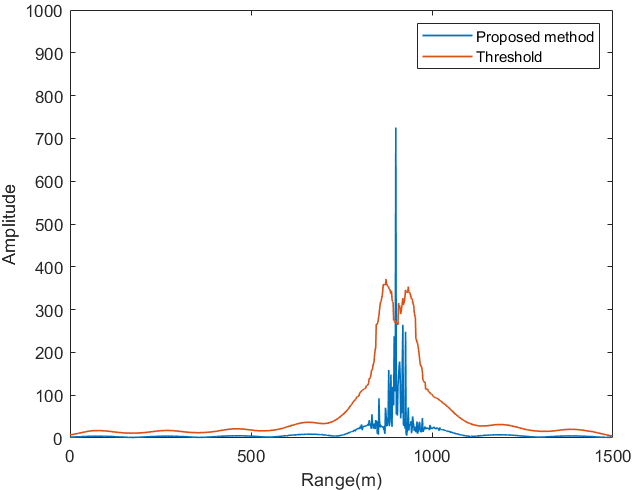}
        \caption{}
    \end{subfigure}
    \caption{Results of CFAR detection. (a) Pulse compression. (b) Energy function. (c) Max-TF. (d) Proposed method.}
    \label{fig11}
\end{figure}

The GLWD image is subjected to M-LLSD, and the result is shown in Fig. ~\ref{fig9}. It can be seen that after detection, the target signal is effectively detected. According to the detection results, a filter is constructed, and filtering is performed on the STFT to suppress jamming. The PC results after jamming suppression is shown in Fig. ~\ref{fig10}. Fig. ~\ref{fig10}(a) and Fig. ~\ref{fig10}(b) are the results of the two comparison methods. It can be seen that the energy function method cannot calculate the low-energy target area because the time domain signal is submerged by noise, which will cause target energy loss and jamming cannot be suppressed cleanly during filtering. The max-TF method mistakenly treats the target signal as jamming during filtering because the STFT image is blurred and the energy of the target and the jamming is the same, resulting in the target being filtered out. The proposed method can effectively remove false targets and noise floors, and the target energy is only slightly lost by about 1dB, which does not affect subsequent target detection, as shown in Fig. ~\ref{fig11}.

\begin{figure}[htbp]
    \centering
    \begin{subfigure}[b]{0.32\textwidth}  
        \centering
        \includegraphics[width=\linewidth]{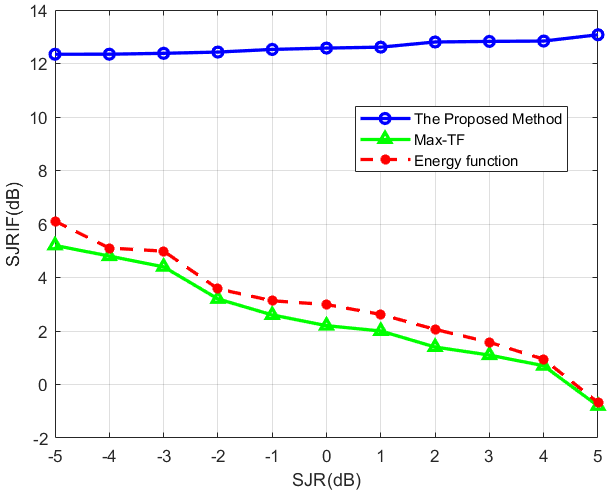}  
        \caption{}
    \end{subfigure}
    \begin{subfigure}[b]{0.32\textwidth}  
        \centering
        \includegraphics[width=\linewidth]{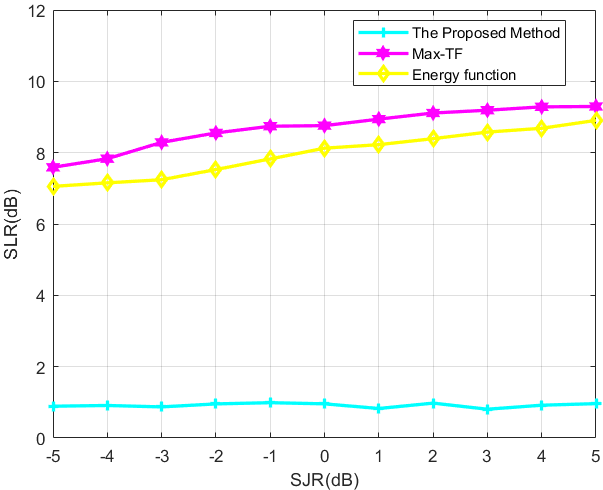}  
        \caption{}
    \end{subfigure}
    \begin{subfigure}[b]{0.32\textwidth}  
        \centering
        \includegraphics[width=\linewidth]{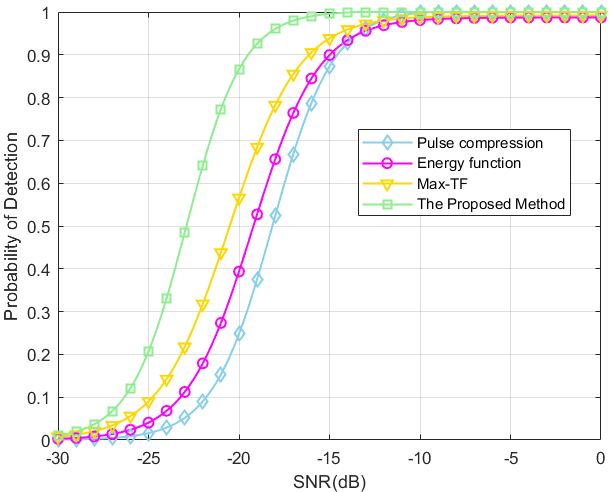}
        \caption{}
    \end{subfigure}
       \caption{Performance curves of the methods. (a) SJRIF curves. (b) SLR curves. (c)  Detection probability curves.}
    \label{fig12}
\end{figure}

Fig. ~\ref{fig12}(a) and Fig. ~\ref{fig12}(b) are the SJRIF and SLR performance curves of different methods, respectively. It can be seen that the SJRIF of the proposed method reaches about 13dB under different SJR, which is significantly higher than that of the comparison method, and the SLR is only about 1dB, while the SJRIF of the comparison method gradually decreases with the increase of SJR, and the SLR is larger.
The target detection performance under different jamming suppression methods is shown in Fig. ~\ref{fig12}(c). The performance of the proposed method is close to the target detection performance without jamming, and has better target detection performance than the comparison method. Compared with the methods in \cite{RN129} and \cite{RN33}, the required SNR is reduced by 4.5 dB and 3.6 dB respectively when the detection probability is about 90\%.
\begin{table}[ht]
\small 
    \begin{center} 
        \begin{minipage}{0.7\textwidth} 
            \centering
            \caption{Simulation Parameters.} 
            \label{table2}
            \small 
            \renewcommand{\arraystretch}{1}  
            \begin{tabular}{@{}p{0.6\textwidth}cl@{}} 
                \toprule
                \text{Parameter} & \text{Value} \\ 
                \midrule
                Waveform pulse duration & 1 $\mu$s \\
Waveform bandwidth & 100 MHz \\
Carrier frequency & 10 GHz \\
Pulse repetition frequency & 5 KHz \\
Targets' location & 900 m \\
Main false targets' locations & 867 m, 927 m \\
ISRJ sampling period & 0.1875$\mu$s, 0.1563$\mu$s \\
ISRJ slice width & 0.125$\mu$s, 0.0625$\mu$s \\
SJR & 4.2 dB, 0 dB \\
SNR & -12 dB \\
                \bottomrule
            \end{tabular}
        \end{minipage}
    \end{center}
\end{table}
\subsection{Compound ISRJ}
 \begin{figure}[htbp]
    \centering
    \begin{subfigure}[b]{0.3\textwidth}  
        \centering
        \includegraphics[width=\linewidth]{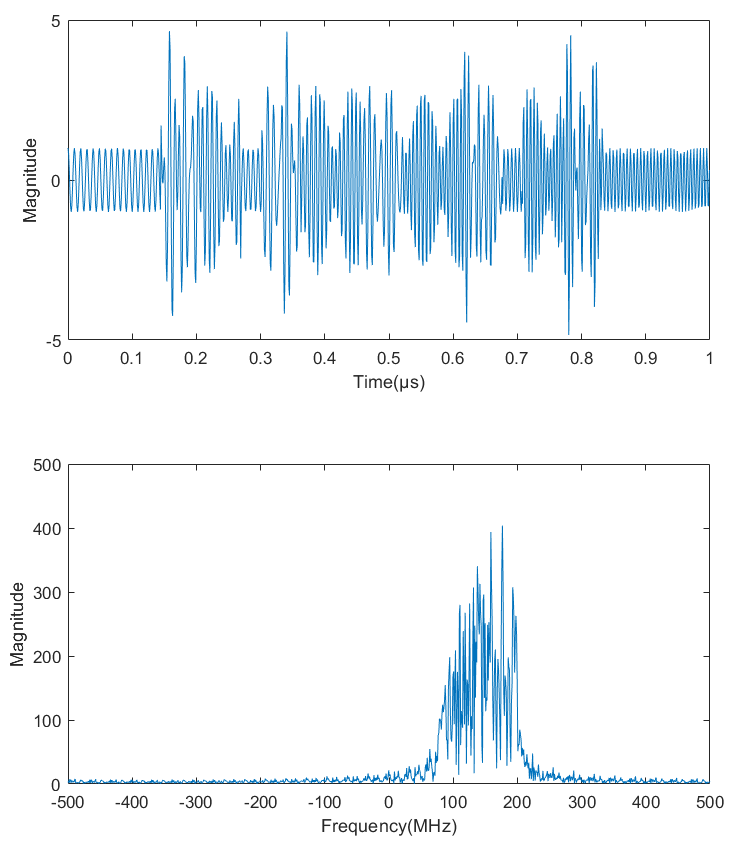}  
        \caption{}
    \end{subfigure}
    \begin{subfigure}[b]{0.3\textwidth}  
        \centering
        \includegraphics[width=\linewidth]{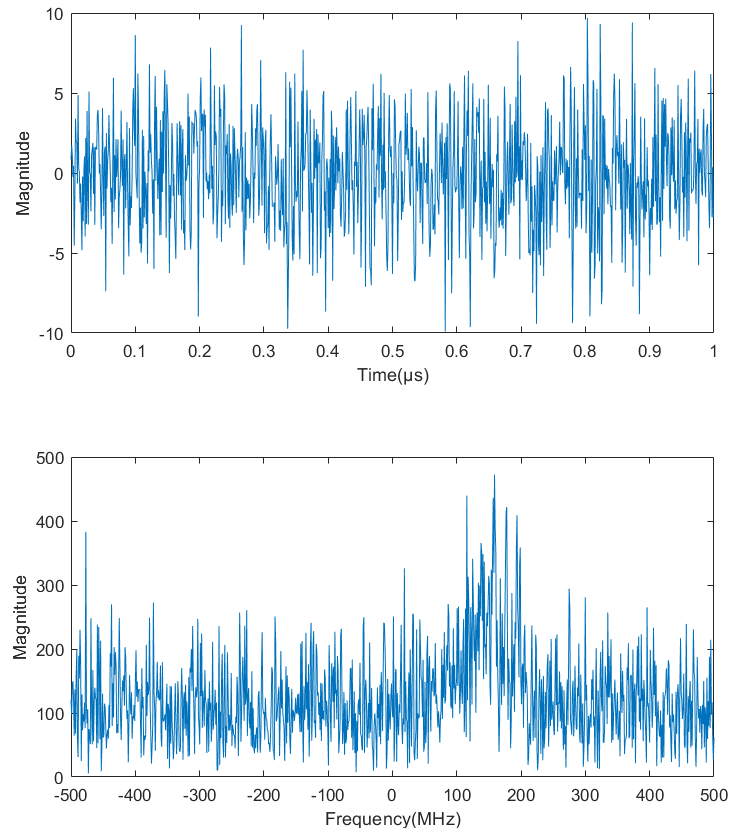}
        \caption{}
    \end{subfigure}
       \caption{Time and frequency domains of signals in compound ISRJ scenario. (a) Target+ISRJ. (b) Target+ISRJ+Noise.}
    \label{fig13}
\end{figure}
\begin{figure}[htbp]
    \centering
    \begin{subfigure}[b]{0.3\textwidth}  
        \centering
        \includegraphics[width=\linewidth]{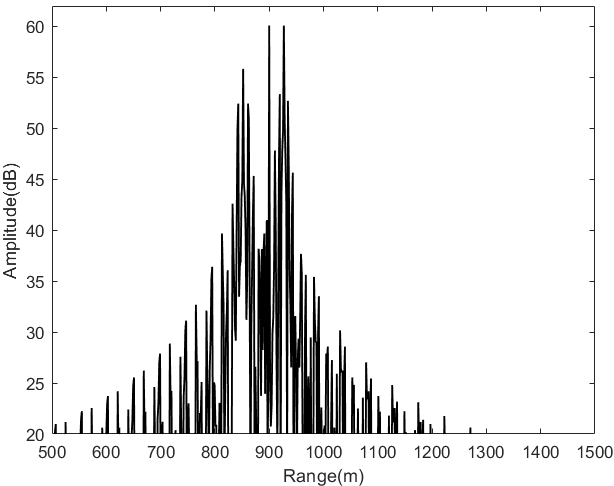}  
        \caption{}
    \end{subfigure}
    \begin{subfigure}[b]{0.3\textwidth}  
        \centering
        \includegraphics[width=\linewidth]{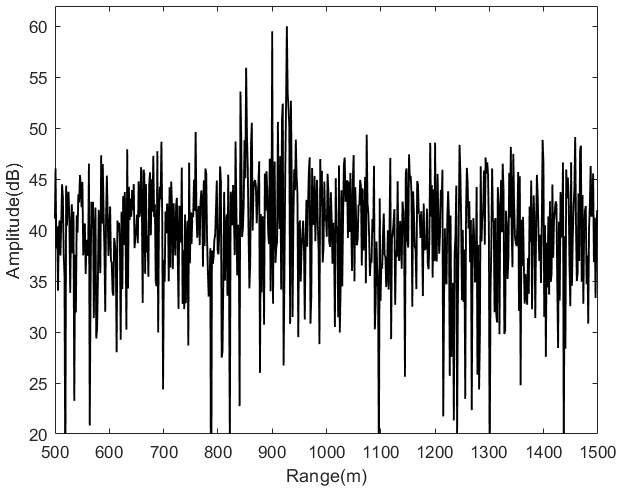}
        \caption{}
    \end{subfigure}
   \caption{Pulse compression in compound ISRJ scenarios. (a) Target+ISRJ. (b) Target+ISRJ+Noise.}
    \label{fig14}
\end{figure}

The difference between this experiment and that in Section \ref{sec:single_isrj} is that an ISRJ is added in front of the target. The relevant parameter settings are shown in Table\ref{table2}. The time domain and frequency domain of the signal are shown in Fig. ~\ref{fig13}, and the PC is shown in Fig. ~\ref{fig14}.

\begin{figure}[htbp]
    \centering
    \begin{subfigure}[b]{0.3\textwidth}  
        \centering
        \includegraphics[width=\linewidth]{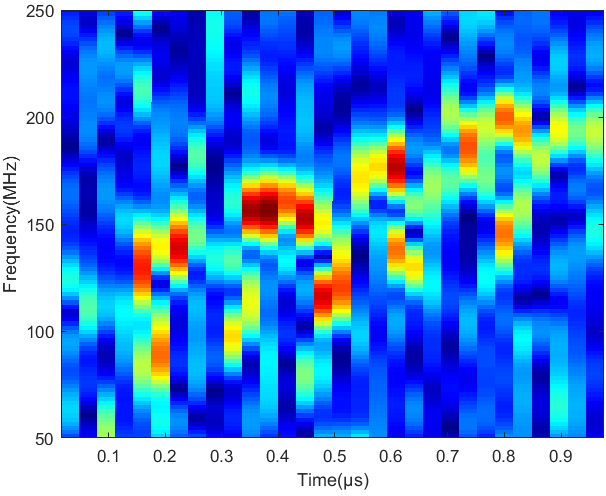}  
        \caption{}
    \end{subfigure}
    \begin{subfigure}[b]{0.3\textwidth}  
        \centering
        \includegraphics[width=\linewidth]{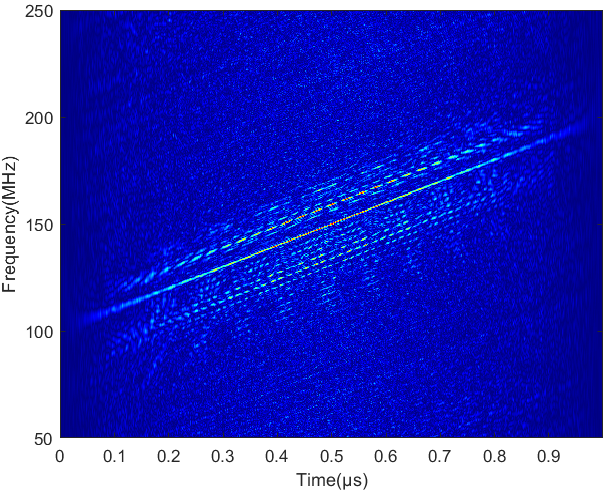}
        \caption{}
    \end{subfigure}
       \caption{Time-frequency image. (a) STFT. (b) GLWD.}
       \label{fig15}
\end{figure}

\begin{figure}
\centerline{\includegraphics[width=0.6\textwidth]{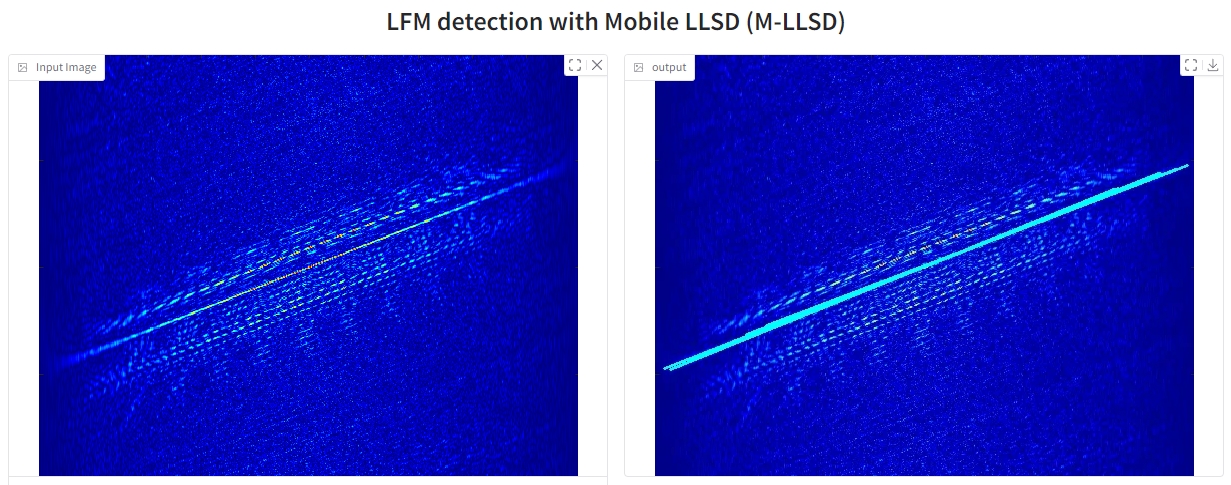}}
\caption{Line detection of target signal.}
\label{fig16}
\end{figure}

\begin{figure}[htbp]
    \centering
    \begin{subfigure}[b]{0.32\textwidth}  
        \centering
        \includegraphics[width=\linewidth]{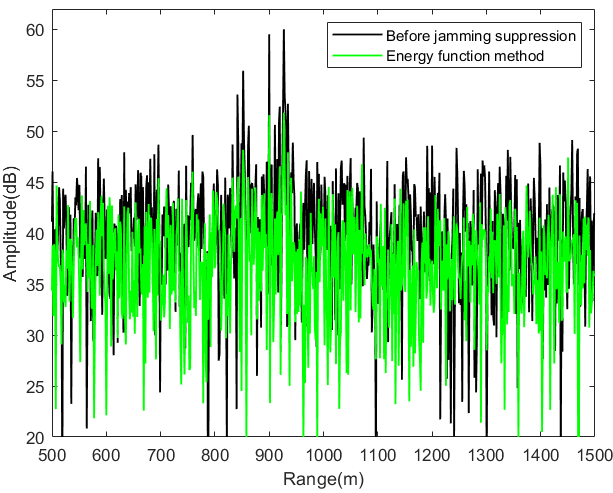}  
        \caption{}
    \end{subfigure}
    \begin{subfigure}[b]{0.32\textwidth}  
        \centering
        \includegraphics[width=\linewidth]{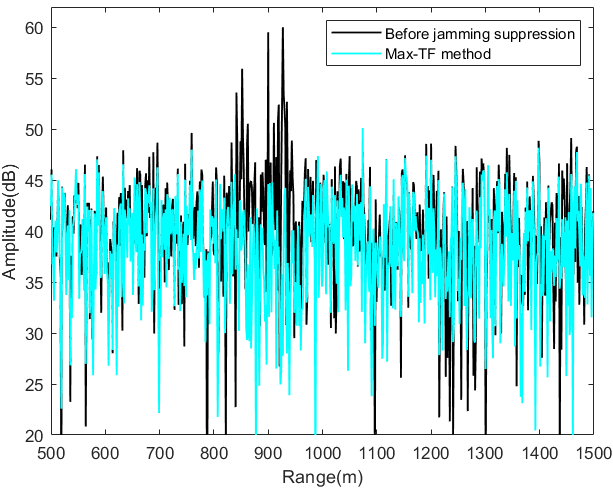}
        \caption{}
    \end{subfigure}
    \begin{subfigure}[b]{0.32\textwidth}  
        \centering
        \includegraphics[width=\linewidth]{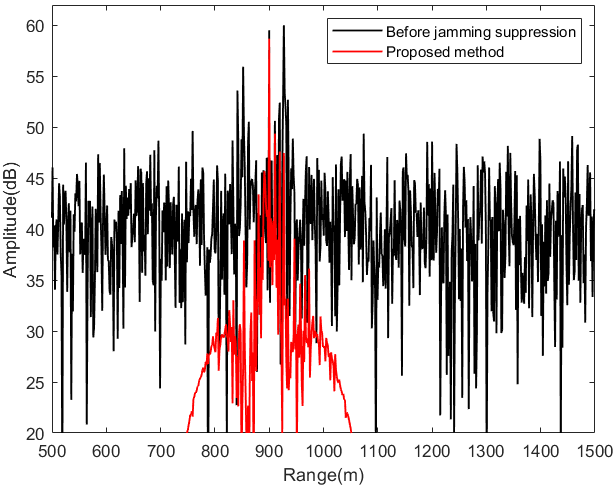}
        \caption{}
    \end{subfigure}
    \caption{Jamming suppression results. (a) Energy function. (b) Max-TF. (c) Proposed method.}
    \label{fig17}
\end{figure}

\begin{figure}[htbp]
    \centering
    \begin{subfigure}[b]{0.24\textwidth}  
        \centering
        \includegraphics[width=\linewidth]{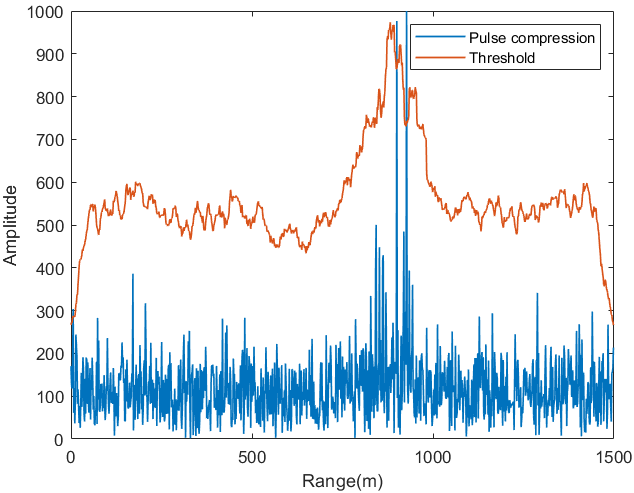}  
        \caption{}
    \end{subfigure}
    \begin{subfigure}[b]{0.24\textwidth}  
        \centering
        \includegraphics[width=\linewidth]{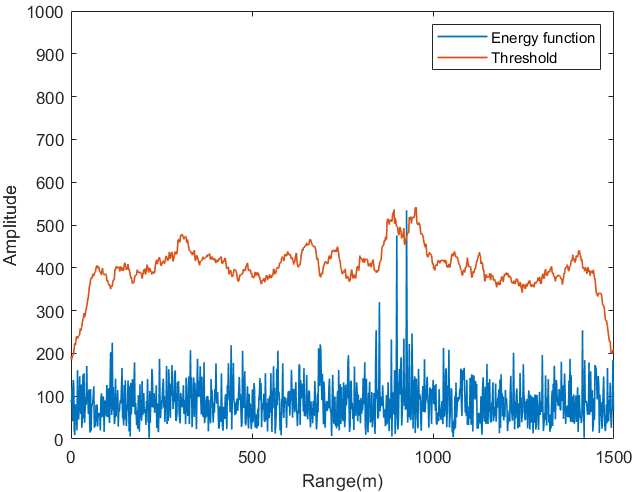}
        \caption{}
    \end{subfigure}
    \begin{subfigure}[b]{0.24\textwidth}  
        \centering
        \includegraphics[width=\linewidth]{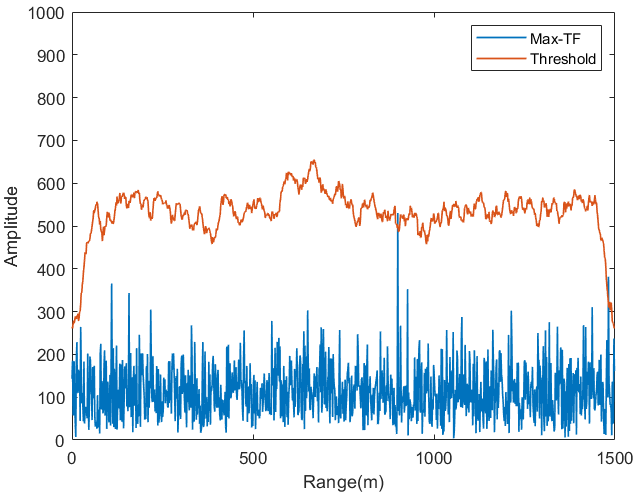}
        \caption{}
    \end{subfigure}
    \begin{subfigure}[b]{0.24\textwidth}  
        \centering
        \includegraphics[width=\linewidth]{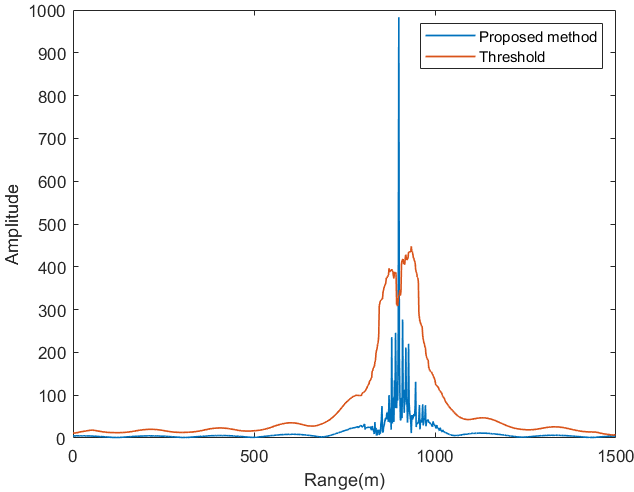}
        \caption{}
    \end{subfigure}
    \caption{Results of CFAR detection. (a) Pulse compression. (b) Energy function. (c) Max-TF. (d) Proposed method.}
    \label{fig18}
\end{figure}

Fig. ~\ref{fig15} is the TF image of the received signal. It can be seen that in this case, the signal aliasing in the STFT domain is serious and the image is blurred, while the GLWD is relatively clear and the signals are clearly separated. The GLWD image is subjected to straight line detection, and the detection result is shown in Fig. ~\ref{fig16}. The target signal is effectively detected, and the TF position is determined, so that a filter is constructed for filtering. The final PC results of the jamming suppression output is shown in Fig. ~\ref{fig17}.
The proposed method can effectively suppress false targets and noise floor, and the target is clearly visible and can be effectively detected, as shown in Fig. ~\ref{fig18}.

\begin{figure}[htbp]
    \centering
    \begin{subfigure}[b]{0.32\textwidth}  
        \centering
        \includegraphics[width=\linewidth]{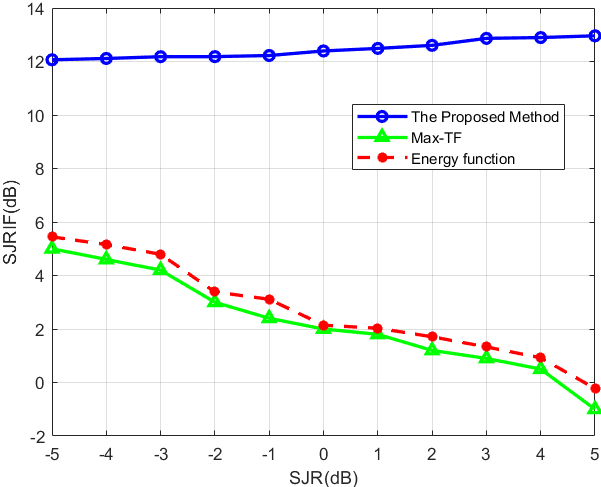}  
        \caption{}
    \end{subfigure}
        \begin{subfigure}[b]{0.32\textwidth}  
        \centering
        \includegraphics[width=\linewidth]{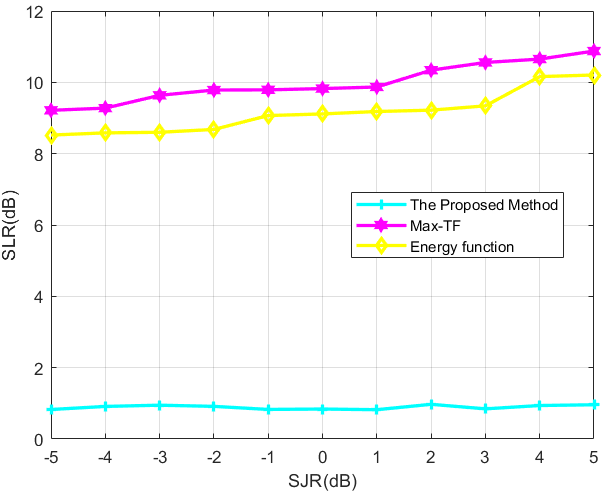}  
        \caption{}
    \end{subfigure}
        \begin{subfigure}[b]{0.32\textwidth}  
        \centering
        \includegraphics[width=\linewidth]{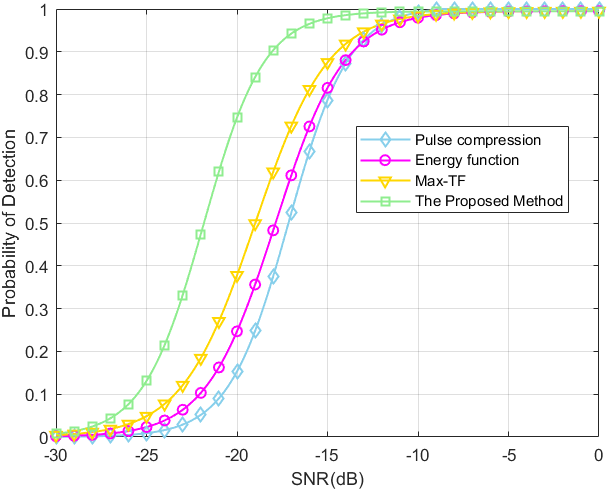}
        \caption{}
    \end{subfigure}
       \caption{Performance curves of the methods. (a) SJRIF curves. (b) SLR curves. (c)  Detection probability curves.}
    \label{fig19}
\end{figure}
Fig. ~\ref{fig19} is the performance curves of different methods. The proposed method performs better than the comparison method in SJRIF and SLR. Fig. ~\ref{fig19}(c) is the target detection performance under different jamming suppression methods. Compared with single ISRJ, the target detection performance is slightly reduced, which is caused by the increase in jamming false targets. However, compared with the methods in \cite{RN129} and \cite{RN33}, the proposed method still has obvious advantages.

\section{Conclusion} \label{sec5}
Existing methods have difficulty in suppressing ISRJ false targets whose energy-frequency is very close to the real target, and the performance will be significantly reduced under low SNR. To address this challenge, this paper proposes an ISRJ suppression method based on linear canonical WD and improved lightweight line detection. This method uses GLWD to process the received signal to improve the TF separation and energy concentration of the signal. By applying M-LLSD to perform line detection on the TF image, the TF position of the target signal can be effectively confirmed. Considering that nonlinear transformations such as WD do not have a direct inverse transformation, a TF filter is constructed based on the mapping relationship between GLWD and STFT, filtered in the STFT domain, and then inversely transformed to obtain the time domain signal to complete jamming suppression. Simulation and experimental results verify the effectiveness of this method in complex ISRJ environments, and it has good robustness and real-time performance.
In summary, this study proposed a new anti-ISRJ method based on linear canonical WD line detection. This algorithm makes up for the shortcomings of existing methods in suppressing ISRJ with energy and frequency similar to the real target under low SNR, fully combines the advantages of traditional signal processing methods and new lightweight intelligent methods, greatly improves the upper limit of TF domain anti-jamming methods, has good practical application prospects, and lays the foundation for further exploration of LCT domain anti-jamming methods.


\section*{Declaration of competing interest}
The authors declare that they have no known competing financial interests or personal relationships that could have appeared to influence the work reported in this paper.

\section*{Acknowledgments}
This work was supported by grants from the National Natural Science Foundation of China [No. 62171041]; and the Natural Science Foundation of Beijing Municipality [No. 4242011].


\bibliography{mybibfile}
\newpage
\end{document}